\documentclass[12 pt]{article}
\usepackage{fullpage}
\usepackage[utf8]{inputenc}
\usepackage{graphicx}
\usepackage{graphics}
\usepackage{float}
\usepackage{color}
\usepackage{mathtools}
\usepackage{amsmath}
\usepackage{amssymb}
\usepackage{etoolbox}
\usepackage{verbatim}
\newcommand{\kstep}[1]{#1_{k+\frac{1}{2}}}
\newcommand{\kstem}[1]{#1_{k-\frac{1}{2}}}
\newcommand{\norm}[1]{|#1|}

\begin{document}

      \begin{titlepage}
      	\begin{center}
     
       	\newcommand{\HRule}{\rule{\linewidth}{0.5mm}} 
       	
       	\center 

       	\textsc{\LARGE University of Science And Technology, Zewail City}\\[1.5cm]
       	
       	\textsc{\Large Aerospace Engineering Department}\\[0.5cm] 
       	\HRule\\[0.4cm]
       	
       	{\Large\bfseries A Geometric Approach to Modeling, Simulation and Control }\\[0.5cm] 
       	{\large  With Application to Control of Rigid Body Dynamical Systems}\\[0.4cm]
       	\HRule\\[1.2cm]
       	
       			\large
       			\textit{Authors}\\
       			\textsc{Mahmoud Abdelgalil \\ Asmaa Eldesoukey \\ Esraa Elshabrawey} 
      
\HRule\\[1.2cm]
       			\large
       			\textit{Supervisor}\\
       			\textsc{Dr. Moustafa Abdallah} %
       			\vfill	
\begin{figure}[H]
	\centering
	\includegraphics[scale =.7]{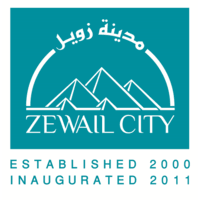}
	\centering
\end{figure}
\center 
\textsc{Senior Design Project submitted for the degree of Aerospace Engineering}
        	
    \end{center}
\end{titlepage}
 \newpage
 \vfill
 \textit{To those who love mathematics and mechanics.}
\newpage

\begin{large}
\textbf{Dedication} \\

\end{large}
Heartfelt thanks go to our supervisor \textbf{Dr.Moustafa Abdallah}, for his guidance and support throughout the work. He spared no time and effort to teach us. This thesis would have been impossible without his continuous aid and support. It has been such an honor to be exposed to this much of academic guidance, insight and genius. Moreover, he provided a healthy work environment through competition and friendship among the team members, as well.  \\

We are profoundly grateful to \textbf{Dr.Haithem Taha} for his encouragements to our research interests. He was the first to inspire us through his work in geometric mechanics and control. Our interactions with him motivated us to pursue our academic path with enthusiasm towards control systems. \\

To the memory of \textbf{Dr.Ahmed Zewail}, we dedicate this work. He established our university in hopes of enhancing higher education in Egypt. His efforts surely paid off and we wished he would be there witnessing the graduation of the first engineering batch of the university.  \\

To our \textbf{families and friends}, we can not be more thankful for their invaluable reassuring  and immersing us with comfort and confidence especially in hard times.

\begin{flushright}
	Asmaa Eldesoukey \\ Mahmoud Abdelgalil \\ Esraa Elshabrawey
\end{flushright}

\newpage
\thispagestyle{plain}
\begin{center}
	\large
	\textbf{A Geometric Approach to Modeling, Simulation and Control}

	\vspace{0.4cm}
	Mahmoud Abdelgalil, Asmaa Eldesoukey, Esraa Elshabrawy
	
	\vspace{0.9cm}
	\textbf{Abstract}\\
	In this work, we utilize discrete geometric mechanics to derive a 2nd-order variational integrator so as to simulate rigid body dynamics. The developed integrator is to simulate the motion of a free rigid body and a quad-rotor. We demonstrate the effectiveness of the simulator and its accuracy in long term integration of mechanical systems without energy damping. Furthermore, this work deals with the geometric nonlinear control problem for rigid bodies where backstepping controller is designed for full tracking of position and orientation. The attitude dynamics and control are defined on \textbf{SO(3)} to avoid singularities associated with Euler angles or ambiguities accompanying quaternion representation.  The controller is shown to track large rotation attitude signals close to $180^\circ$  achieving almost globally asymptotic stability for rotations. A Quad-rotor is presented as an example of an under-actuated system with nonlinear model on which we apply the backstepping control law. In addition, an aerodynamic model aiming at deriving the aerodynamic forces and torques acting on rotors is added for realistic simulation purposes and to testify the effectiveness of the derived control method.
\end{center}
\newpage
\tableofcontents
\newpage
	\section{Motivation And Literature Survey}
	\subsection{Discrete Variational Mechanics}
	The main advantage of using the formalism of discrete variational mechanics is seeking more accurate solutions to non-linear differential equations that are derived from variational principles, which are essentially all mechanical systems. In \cite{marsden2001}, the authors show a graph of a simulation of a mechanical system and compare the performance of a non-variational integration scheme (Runge-Kutta) and a variational integrator scheme (Newmark). The system has no dissipative forces, which means that the total energy of the system should remain constant over time. 
	\begin{figure}[H]
		\centering
		\includegraphics[width=0.5\textwidth]{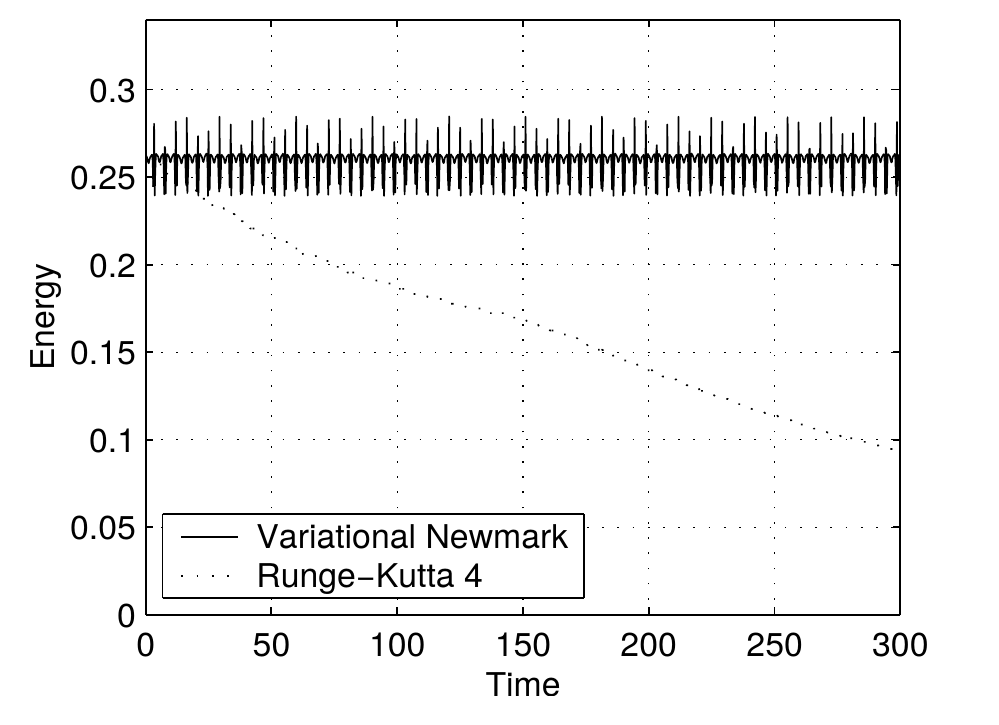}
		\caption{Comparison between Runge-Kutta and Newmark integration schemes}
	\end{figure}
	The striking aspect of the graph is that normal integration using 4th order runge-kutta scheme is shown to decay the energy dramatically when integrated for a long time. This may for example pose a problem in the case of simulation of a control system after designing the control. The system may not include damping, and yet due to the presence of this artificial damping introduced by the integration scheme, the control design engineer might be satisfied with the performance shown by the simulator thinking that the system will behave in the same way, while in fact the presence of this artificial damping generally introduces some artificial stability in the system and thus the actual performance of the control law will be worse when applied to the real system.  \\
	\subsection{Geometric Control}
	Most of the mechanical systems have configuration spaces of smooth non-Euclidean manifolds such as  Lie groups. Generally this applies to dynamical systems with rotational DOF. For instance,
	consider the case of a two link manipulator (double pendulum) with two degrees of freedom $\theta_1, \theta_2$,
	where $\theta_1, \theta_2 \in S^1$. Hence,
	the configuration space is a torus $T: S^1 \times S^1$ and not $\mathbb{R}^2$. 
	\begin{figure}[H]
		\centering
		\includegraphics[scale =.8]{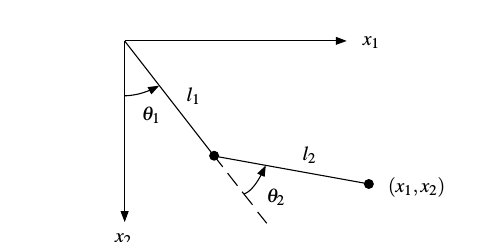}
		\caption{two link manipulator, $\theta_1, \theta_2 \in S^1$ \cite{schaft}.}
		\centering
	\end{figure}
	Another example is a rigid body rotation evolving on a smooth manifold which is the Lie group of the $3 \times 3$ orthonormal matrices: $\textbf{SO(3)}$. \\
	Such systems and many other mechanical ones need to be treated with a special class of controllers for their strong nonlinearity. Geometric control is the class of control that deals with systems whose configuration spaces are smooth manifolds utilizing the language of differential geometry. \\
	
	The conventional way to handle system coordinates evolving on manifolds is to \textit{locally} mark the change of these configuration variables. This legitimizes considering the system coordinates Euclidean ones. However, addressing the control problem this way loses the essential properties assigning the global nature of systems. An unveiled fact in \cite{Bhat} is that there is no globally asymptotically stable equilibrium  of a continuous dynamical system on a compact manifold with a continuous control vector field. Owing to the fact that it is not tenable to parametrize a whole manifold with a single chart. Another interesting result is that controllers acquired from systems parametrized with Euclidean coordinates upwind the system dynamics in case of deflections such as $ \pm 2 \pi, \pm 4\pi ,\pm 6 \pi..$ rotations though the system rests actually at its equilibrium (origin). \\
	
	Another merit of studying geometric control theory is to possess tools of nonlinear controllability in order to perform forbidden motions, unactuated, with interactions between admissible control vector fields. For instance, the work of Hassan and Taha in \cite{Hassan} revealed new rolling/yawing mechanisms using certain combination of aileron and elevator deflections. Similarly, a new pitching mechanism is acquired with certain interaction between elevator motion and throttle deflections. Moreover, Burcham et al. in \cite{Pappas} presented the capability of controlling a spacecraft with failure in a couple of gas jets. Hereby, the geometric control tools assist in the fault tolerance analysis associated with actuation failure.\\
	
	Such interesting findings motivated this dissertation. The derived backstepping geometric controller is expressed in \textbf{SO(3)} manifold and is shown to track large rotation attitude signals and desired position ones achieving almost globally asymptotic stability. \\
	
	In \cite{bouabdallah}, Bouabdallah and Siegward suggested two types of nonlinear controllers: backstepping and sliding mode for Quad-rotor tracking. The suggested controllers did not represent the attitude dynamical model nor the controller in \textbf{SO(3)} yet rotation is parametrized using Euler angles. The controller showed acceptable performance for stabilization in case of backstepping technique for angles perturbations close to $45^\circ$. 
	In \cite{Lee}, Lee developed a geometric controller defined globally on \textbf{SO(3)}. The proposed rotation error function  utilizes the Riemannian metric consistent with the rotation group.  For large initial attitude errors, the controller fulfilled efficiently tracking the desired attitudes with small steady state errors.
	In \cite{Raj}, Raj et al. applied the geometric backstepping control theory on a small-scale rotary wing aircraft for tracking. The controller achieved desirable performance especially with roll aggressive (high frequency sinusoidal signal, very large initial error and high roll rate) maneuver with acceptable control input.

	\newpage
		\section{Background}
		In preparation for later chapters, the reader is recommended to read this chapter as it contains basic however sufficient background in differential geometry and algebra to proceed with the geometric mechanics and control to come. 
		\subsection{Manifolds}
		Crudely, manifolds are abstract spaces that locally look like linear spaces. If the manifold is differentiable/differential, it looks so similar to linear spaces that we can apply the principles of calculus and other operations like addition as we know it in linear spaces. A key concept of defining a manifold is mapping to $\mathbb{R}^n$ via coordinate charts. 
		\subsubsection{Coordinate Charts}
		Let $M$ be a manifold and $U$ be a set on $M$, a coordinate chart is  the set $U$ along with a homeomorphic (continuous with continuous inverse) map $\phi: U \longrightarrow \phi(U) \subset \mathbb{R}^n$ where $\phi(U)$ is an open set in $\mathbb{R}^n$. Two charts $(U_1,\phi_1), (U_2,\phi_2)$ where $U_1, U_2$ are not disjoint are said to be compatible if $\phi_1(U_1 \bigcap U_2) , \phi_2(U_2 \bigcap U_1)$ are open sets in $\mathbb{R}^n$ and their compositions $\phi_1 \circ \phi_2^{-1}, \phi_2 \circ \phi_1^{-1}$ are $C^\infty$. 
		\begin{figure}[H]
			\centering
			\includegraphics[scale =.8]{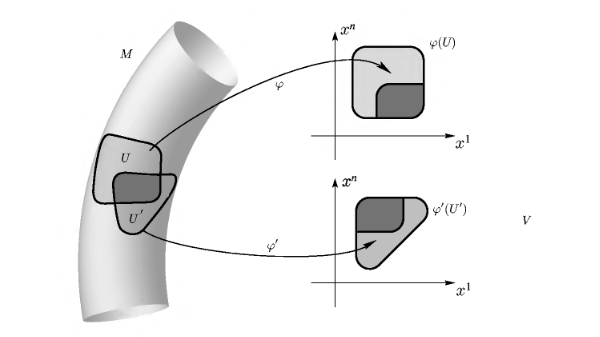}
			\caption{A manifold with coordinate charts, \cite{Marsden}.}
			\centering
		\end{figure}
		A differentiable manifold is now can be defined as space that \textit{everywhere} covered with a collection of compatible charts (Atlas). The maximal atlas, the atlas contains all the compatible coordinate charts, defines the differentiable structure associated with the differential manifold. We can say the differentiable manifold is diffeomorphic to $\mathbb{R}^n$ where a diffeomorphism is a map from the differentiable manifold $M$ to $\mathbb{R}^n$ with an inverse such that the map and inverse are smooth.
		\subsubsection{Curves}
		Let $M$ be a manifold, a curve on $M$ is a map $c: t \subset \mathbb{R}  \longmapsto M$. A tangent vector of a curve $c$ is the velocity at a point $p$.
		\begin{figure}[H]
			\centering
			\includegraphics[scale =.7]{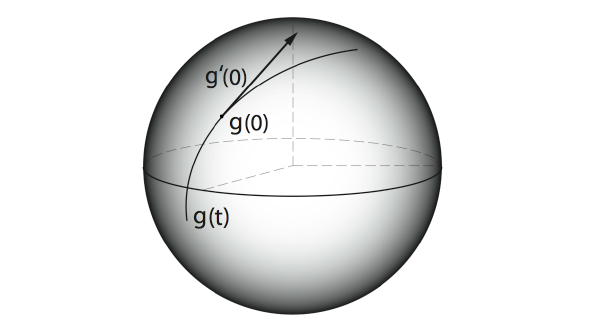}
			\caption{$S^2$ manifold with curve $g(t)$ and a tangent vector at $g(0)$, \cite{holm}.}
			\centering
		\end{figure}
		\subsubsection{Tangent and Cotangent Spaces and Bundles}
		Let $M$ be a smooth/differentiable manifold, the tangent space of $M$ at any point $p$ denoted by $T_p M$ is the space of all tangent vectors through $p$. Any tangent space at a point is a vector space isomorphic to $\mathbb{R}^n$, its map to $\mathbb{R}^n$ is bijective and preserves group operation.\\
		
		The tangent bundle of $M$ denoted $TM$ is the disjoint union of all the tangent spaces for all points $p \in M$. For n-dimensional manifold, the tangent bundle has a dimension of $2 n$. \\
		The tangent bundle projection is a map $\lambda: TM \longrightarrow M$ yielding the base point $p$ corresponding to each tangent vector. \\
		
		The tangent space  $T_p M$ of $M$ at any point $p$  has a dual space denoted by $(T_p M)^*$ called the cotangent space.
		The cotangent bundle of $M$ denoted $(TM)^*$ is the union/vector bundle of all the cotangent spaces. The dual of a vector space is the space of all linear functionals along with the vector space.  \\
		\begin{figure}[H]
			\centering
			\includegraphics[scale =.7]{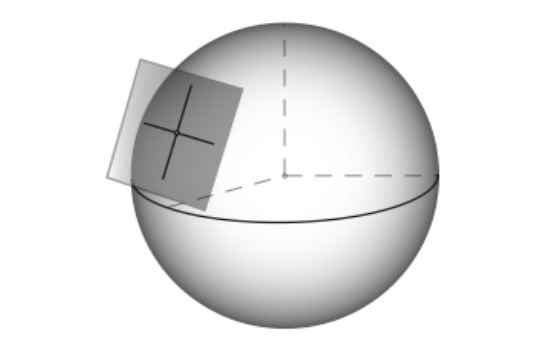}
			\caption{Tangent space $T_p S^2$ $\mathbb{R}^2$ to $S^2$ manifold, \cite{holm}.}
			\centering
		\end{figure}
		\subsection{Directional and Lie derivatives}
		\subsubsection{Vector Fields}
		A vector field $X$ of a Manifold $M$ is an assignment of a tangent vector at each point on the manifold. If the manifold has coordinates $(x_1,x_2,x_3,..., x_n)$ a vector field is defined as $X(x_1,x_2,x_3,..., x_n)$.
		\subsubsection{Flow Maps}
		For a vector field $X$ of a Manifold $M$, $\phi_t^X(p): M\longrightarrow M$ is a flow map which measures the transition of a point $p$ over a curve $C(t)$ for some time interval $(t_i,t_f)$.
		\subsubsection{Directional Derivative}
		Let $f$ be a differentiable multi-variable function, and a vector $v \in T_p M$. The directional derivative $v(f): C^\infty(M) \longrightarrow \mathbb{R}$ is the change of this function in the direction of the vector $v$.
		\begin{eqnarray}
		v(f(X)) &=& v_i \frac{\partial f(X)}{\partial X} \\
		v(f(X)) &=& v. \nabla f
		\end{eqnarray}
		where $\nabla f$ is the gradient of $f$.
		\subsubsection{Lie Derivative}
		If the change in a multi-variable differentiable function $f$ is evaluated along a vector field $X(x)$, we call it the Lie derivative denoted $L_f X:C^\infty(M)\longrightarrow C^\infty(M) $.
		\begin{eqnarray}
		L_f &=&  \sum_{i}^{n} X_i \cdot\frac{\partial f }{\partial x_i} \\
		L_f &=&  X.\nabla f
		\end{eqnarray}
		where $x_i$ are the coordinates of $X$.
		\subsubsection{Lie Algebra}
		A Lie algebra defines the set of all smooth vector fields on $M$. A vector space $V$ with binary map $V \times V \longrightarrow V$ defines a Lie algebra if it satisfies the properties of bi-linearity, skew symmetry and Jacobi identity.
		\subsubsection{Lie Bracket}
		The Lie derivative of a vector field $X_1$ w.r.t another vector field $X_2$ is called the Lie bracket $[X_1,X_2]: \Gamma(M) \times \Gamma(M) \longrightarrow \Gamma(M)$, where $\Gamma(M)$ is the set of all vector fields on $M$.
		\begin{eqnarray}
		[X_1,X_2] &=& \sum_{j}^{n}(\sum_{i}^{n} \frac{\partial X_{2,j}}{\partial x_i} X_{1,i}-\frac{\partial X_{1,j}}{\partial x_i} X_{2,i}) \frac{\partial}{\partial x_j}
		\end{eqnarray}
		where
		Equivalently,
		\begin{eqnarray}
		[X_1,X_2] &=& L_{X_1} L_{X_2} - L_{X_2} L_{X_1} 
		\end{eqnarray}
		For any vector fields $X_1,X_2,X_3$ on $M$, the following three statements are true. 
		\begin{align}
		[X_1,X_2] &= -[X_2,X_1] \ \ \ \ \ \  &\text{Skew-Symmetry} \\
		[X_1+X_2,X_3] &= [X_1,X_3]+[X_2,X_3]  \ \ \ \ \ \ \  &\text{Linearity} \\
		[X_1,[X_2,X_3]]&+[X_2,[X_3,X_1]]+ [X_3,[X_1,X_2]] =0       \ \ \ \ \ \  &\text{Jacobi Identity}
		\end{align}
		It is proven that 
		\begin{eqnarray}
		t^2 [X,Y](p) &=& \phi^{-Y}_{t} \circ  \phi^{-X}_{t} \circ  \phi^{Y}_{t} \circ  \phi^{X}_{t} (p)
		\end{eqnarray}
		If we start motion at some point $p$ along the vector field $X$ then switching to $Y$ after time t, then back to $-X$ after time step t and again to $Y$, Lie brackets can check if we would get back to $p$ or not;
		$[X,Y]=0$ if and only if $ \phi^{Y}_{t} \circ  \phi^{X}_{t} (p) = \phi^{-Y}_{t} \circ  \phi^{-X}_{t}(p)$. \\
		This has very important implications in the control theory in terms of the importance of actuation order. \\
		Another usage of a Lie bracket is that it signifies the interaction between vector fields in a nonlinear dynamical system. A Lie bracket can measure the possibility of motion in unactuated directions if $[X_1,X_2] \notin span\{X_1,X_2\}$.\\
		Consider the following control affine system model of a car:
		\begin{align}
		\dot{x} &= u_1 g_1(x)+u_2 g_2(x) \\
		\begin{bmatrix} 
		\dot{x} \\ \dot{y} \\ \dot{\theta} \end{bmatrix} &= u_1 \begin{bmatrix}
		\cos \theta \\ \sin \theta \\ 0
		\end{bmatrix} +u_2  \begin{bmatrix}
		0 \\ 0 \\ 1
		\end{bmatrix}
		\end{align}
		where $u_1 = \sqrt{\dot{x}^2+\dot{y}^2}, \ u_2 =\dot{\theta}$ and the configuration space is $\mathbb{R}^2 \times S^1$.
		\begin{eqnarray}
		[g_1,g_2] = \begin{bmatrix} 
		0 & 0 & \sin \theta \\
		0 & 0 &  -\cos \theta \\
		0 & 0 & 0
		 \end{bmatrix} 
		 \begin{bmatrix} 
		 0 \\ 0 \\ 1
		 \end{bmatrix} =  \begin{bmatrix} 
		 \sin \theta  \\-\cos \theta   \\ 0
		 \end{bmatrix}
		\end{eqnarray}
		The result of $[g_1,g_2]$ is a vector field perpendicular to $g_1$ and $g_2$. It shows the motion in the unactuated direction (side motion) perpendicular to the current orientation. The capability of car side motion is  a merit of the nonlinear dynamical systems. With the right combination of admissible inputs this side motion is possible. If the car moves forward then rotates counterclockwise $90^\circ$ then moves backward then rotates clockwise $90^\circ$ and finally moves backward, it is equivalent to move in the side direction.

\subsection{The Rotation Group SO(3)}
		We have seen that the degrees of freedom associated with the relative orientation of body fixed axes with respect to reference axes are described by rotation matrices, hence the name rotational degrees of freedom. In this section we will discuss the rotation group in a more abstract setting where we deal with it as a Lie group. Treating rotation matrices in this setting allows one to perform operations such as differentiation and finding mean (which will prove to be useful in later chapters) in a manner that is consistent with the group structure of \textbf{SO(3)}. Of course, this has the advantage, of being able to treat systems modeled on \textbf{SO(3)} (all systems involving rigid body motion!) in a global manner and be able to evaluate statements as strong as "almost global stability" of a control system designed for a rigid body. 
		we provide a summary of some of the important properties of the rotation matrices, which are elements of the rotation group \textbf{SO(3)}, which turns out to be a Lie group. We begin by noting that, being inner-product preserving, rotation matrices acting on the $\mathbb{R}$ are orthogonal matrices. Particularly, for a rotation matrix \textbf{A}, we have:
		$$ \textbf{A} \textbf{A}^T = \textbf{A}^T \textbf{A} = \mathbb{I}_{3\times 3} $$
		Differentiating the above, we get:
		\begin{eqnarray}
		d\textbf{A} \textbf{A}^T + \textbf{A} d\textbf{A}^T = d\textbf{A}^T \textbf{A} + \textbf{A}^T d\textbf{A} = \textbf{0}_{3\times 3}
		\end{eqnarray}
		Equivalently, we can write:
		\begin{eqnarray}
		d\textbf{A} \ \textbf{A}^T = - \textbf{A} \ d\textbf{A}^T = (d\textbf{A} \textbf{A}^T)^T
		\end{eqnarray}
		Which amounts to saying that the matrix $d\textbf{A} \textbf{A}^T$ is skew symmetric. 
		
		We next show that any proper rotation matrix can be expressed as the exponential of a skew-symmetric matrix. We begin the argument by making the assumption that a matrix $\Theta$ is skew-symmetric and that its exponential is equal to another matrix \textbf{Q}:
		\begin{eqnarray}
		 \textbf{Q} = e^{\Theta} = \mathbb{I}_{3\times 3} + \frac{\Theta}{1!}+ \frac{\Theta}{2!} + \frac{\Theta^3}{3!} + \frac{\Theta^4}{4!} + \; . \ . \ . \ . \ 
		\end{eqnarray}
		What we want to prove is:
		\begin{eqnarray}
		Q Q^T = \mathbb{I}_{3\times 3}
		\end{eqnarray}
		We have:
		\begin{eqnarray}
		 Q^T = (e^\Theta)^T = \mathbb{I}_{3\times 3} + \frac{\Theta}{1!}^T+ \frac{\Theta^2}{2!}^T + \frac{\Theta^3}{3!}^T + \frac{\Theta^4}{4!}^T + \; . \ . \ . \ . \ .
		\end{eqnarray}
		which is in fact equal to:
		\begin{eqnarray}
		Q^T = e^{\Theta^T} = e^{-\Theta}
		\end{eqnarray}
		Thus, we have: 
		\begin{eqnarray}
		Q Q^T = e^{\Theta} e^{-\Theta} = e^{\textbf{0}_{3\times 3}} = \mathbb{I}_{3\times 3}
		\end{eqnarray}
		which is what was to be proven. 
		
		Next we would like to find the derivative of a rotation matrix. We utilize the result just proven and write:
		\begin{eqnarray}
		 d\textbf{A} = d(e^{\Theta})  = d\Theta \ e^{\Theta} = d\Theta A
		\end{eqnarray}
		 	
		The results proven above have deeper meanings when considered in the context of Lie groups. We stated in the beginning of this section that the rotation group \textbf{SO(3)} is a Lie group, which roughly means that its group action can be infinitesimal, i.e the group is in fact a manifold (particularly, it is a connected compact 3-manifold). In this setting, the exponential map proven above is in fact the map that relates the rotation Lie group \textbf{SO(3)} to its associated Lie Algebra $\mathfrak{so}(3)$:
		$$ exp: \mathfrak{so}(3) \rightarrow \textbf{SO(3)} $$
		
		The Lie algebra of the rotation group corresponds to the space of infinitesimal rotations around identity. We use this property of the Lie Algebra to perform differentiation on the manifold corresponding to the Lie group. To see how, we first by considering a rotation matrix that is the exponential of a very small skew-symmetric matrix:
		\begin{eqnarray}
			Q = \mathbb{I}_{3\times 3} + \frac{d\Theta}{1!} + \mathcal{O}(d^2\Theta) 
		\end{eqnarray}
		Now normal differentiation would correspond to 
		\begin{eqnarray}
		 Q - \mathbb{I}_{3\times 3} = d\Theta + \mathcal{O}(d^2\Theta) 
		\end{eqnarray}	
		Hence, we can see that very small (infinitesimal) skew symmetric matrices represent (in the limiting behavior) the difference normal euclidean difference. Using this fact, we can define the derivative of a rotation matrix as we have previously did as:
		$$ dQ = d\Theta \ Q $$
		where $d\Theta$ is an infinitesimal skew-symmetric matrix that amounts to the infinitesimal perturbation of the matrix Q. This infinitesimal perturbation belongs to the Lie algebra of the rotation group, and hence the meaning we wanted to illustrate. 
		\newpage
		\section{Rigid Body Motion}
		\subsection{Degrees of Freedom of a Rigid Body}
		We know from experience that a rigid body has 6 degrees of freedom. In this section we present an argument for  why this is true. Consider a system of N particles moving in 3D space. Define a reference axes in the fixed (inertial) space with unit vectors \textbf{i}, \textbf{j} and \textbf{k}, orthogonal to each other. Thus, to specify the positions of all the particles in the system with respect to this inertial space we need 3\textbf{N} coordinates. This is the most general case for a system of \textbf{N} particles. \\
		
		In the case of a rigid body, the particles in the system are not entirely free to move; they have to move in such a way that keeps the relative distance between each two particles fixed. Denoting the relative distance between the ith and jth particles by $r_{ij}$, the aforementioned constraint is expressed as:
		\begin{equation}\label{eq:cnstr}
		r_{ij} = c_{ij} \quad \forall i \ , j \in \{1,2,3, ..., N\}, \; i\neq j
		\end{equation} 
		where $c_{ij}$ is a constant. Using simple combinatorics, we can see that the number of constraints of the form in  Eq.\ref{eq:cnstr} is equal to $\frac{1}{2} \textbf{N}(\textbf{N}-1)$. This is a far greater number than 3\textbf{N} and hence it is clear that subtracting the number of constraints from the total number of degrees of freedom will not lead us to a correct result. 
		
		In fact, not all constraints of the form in Eq.\ref{eq:cnstr} are independent. To see why, assume we know the position of 3 non-co-linear points that belong to the body. Because the system of particles that constitute a rigid body moves such that the relative distance between particles are constant, knowing the positions of 3 non-collinear points on the body fully specifies the positions of all other particles in the system. Thus, we already know that the number of degrees of freedom of a rigid body is less than or equal to 9, which is the number of coordinates needed to specify the positions of these 3 non-collinear points. However, the positions of these 3 points are themselves not independent. In fact, denoting these points by \textbf{1}, \textbf{2} and \textbf{3} we have:
		\begin{equation}\label{eq:cnstr_red}
		r_{ij} = c_{ij} \quad \forall i \ , j \in \{1,2,3\}, \; i\neq j
		\end{equation} 
		This reduces the degrees of freedom from 9 to 6 which is consistent with our intuition. 
		
		\clearpage
		\subsection{Description of Rigid Body Motion}
		We know from the previous section that a rigid body has 6 degrees of freedom. If we attach a Cartesian set of coordinate axes to a point in the body such that it remains fixed in the body, we can specify the position of any particle in the body with respect to these set of axes. Hence, knowing how to relate positions in these body fixed axes to positions in the reference axes is sufficient to describe the degrees of freedom of the body. Obviously, 3 coordinates are needed to specify the position of the origin of the body fixed axes with respect to the reference axes. Then, the remaining 3 degrees of freedom are associated with the relative orientation of the body fixed axes with respect to the inertial axes. 
		\begin{figure}[H]
			\centering
			\includegraphics[width=0.5\textwidth]{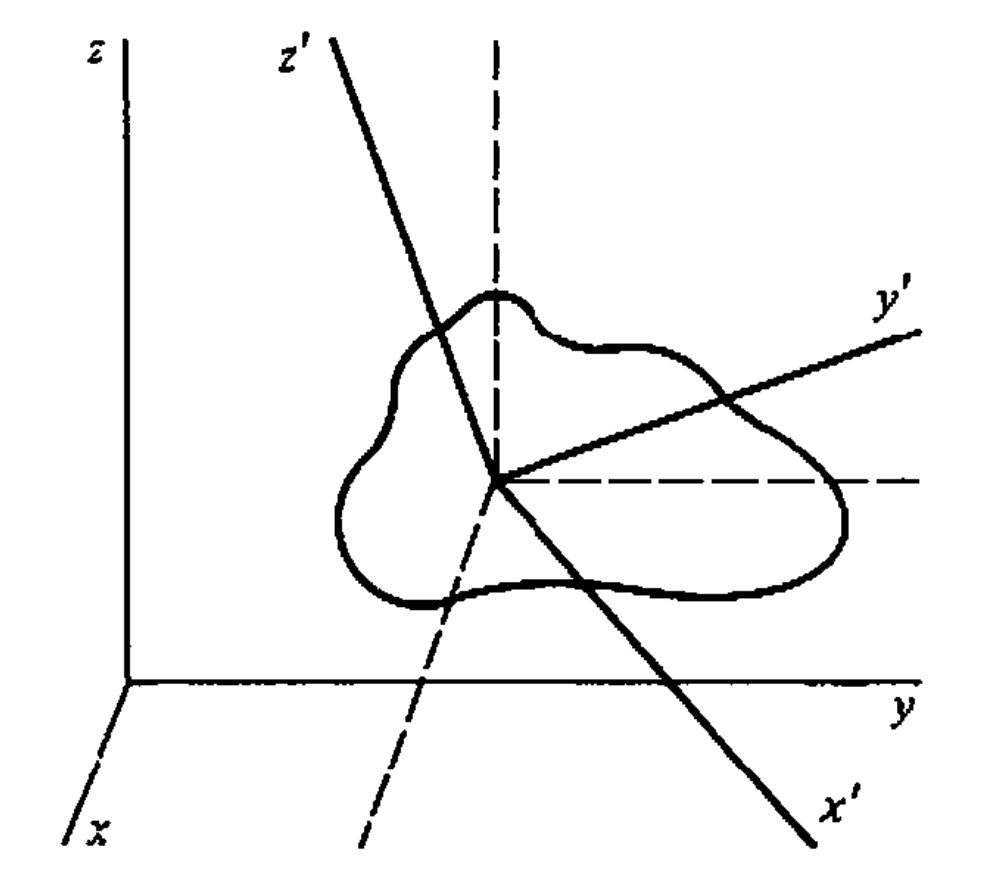}
			\label{fig:body_axes}
		\end{figure}
		It turns out that choosing 3 coordinates to specify the orientation is very tricky. More on this will come later. The above framework is depicted in figure \ref{fig:body_axes}
		
		Since we know that whatever the motion of a rigid body is, it must leave distances between body material points fixed, then we know that the relation between the body axes and reference axes is specified by a rotation matrix. Furthermore, we require the motion to be continuous, which comes from a physical argument, we further constrain ourselves to proper rotation matrices. The group of proper rotation matrices is denoted \textbf{SO(3)} and is known to be a Lie group. \\
		
		So far, we  have reached the conclusion that the number of degrees of freedom of the rigid body is 6 and that the general motion of a rigid body occurs on the space consisting of the Cartesian product of $\mathbb{R}^3$ (to which the position of the origin of the body fixed axes belong) and \textbf{SO(3)} (which is the space of all possible proper rotations that can happen between the body axes and reference axes). Denote this space with \textbf{Q} we can write:
		$$ q(t) \in \textbf{Q}= \mathbb{R}^3 \times \textbf{SO(3)} $$
		where q(t) is the path followed by a rigid body and \textbf{Q} is the configuration space of a rigid body undergoing general motion. 
		\subsection{Equations of Motion}
		In this section, we derive the equations of rotational motion of a rigid body using the Lagrangian formalism. We know that the configuration space of rotational degrees of freedom is the group of rotation matrices \textbf{SO(3)}. Define the Lagrangian as a mapping $\mathcal{L}: \mathfrak{so}(3)\times \textbf{SO(3)} \rightarrow \mathbb{R}$, such that:
		\begin{equation}
		\mathcal{L} = \frac{1}{2}\omega^T J \omega - U(T)
		\end{equation}
		where $\omega$ is the angular velocity of the rigid body expressed in body coordinates, J is the inertia tensor also expressed in body coordinates, $T$ is the rotation matrix that describes the current attitude of the rigid body with respect to an inertial frame and $U$ is the attitude dependent potential energy. Define the action integral as:
		\begin{equation}
		S(\omega,T) = \int_0^T \mathcal{L} dt
		\end{equation}
		From the principle of least action we know that:
		\begin{equation}
		\delta S = \delta \int_0^T \mathcal{L} dt = \int_0^T \left[D_T \mathcal{L}\cdot \delta T + D_{\omega} \mathcal{L}\cdot\delta\omega\right]dt
		\end{equation}
		From the properties of \textbf{SO(3)}, we have:
		\begin{align}
			 \dot{T} &= T \widehat{\omega} \\
			 T^T\dot{T} &= \widehat{\omega} 
		\end{align}
		Thus we can write:
		\begin{eqnarray}
		(\delta T)^T \dot{T} + T^T \delta\dot{T} = (\delta T)^T \dot{T} + T^T \dot{\delta T} 
		\end{eqnarray}
		and we know that:
		\begin{eqnarray}
		\delta T = T \widehat{\delta\theta}
		\end{eqnarray}
		Hence:
		\begin{eqnarray}
			 (\delta T)^T \dot{T} + T^T \delta\dot{T} &=& - \widehat{\delta\theta} T^T\dot{T} + T^T \left(\dot{T}\widehat{\delta\theta} + T^T \widehat{\dot{\delta\theta}}\right) = \widehat{\dot{\delta\theta}} + \widehat{\omega} \widehat{\delta\theta}  - \widehat{\delta\theta} \widehat{\omega} \\
			  \widehat{\delta\omega} &=&  \widehat{\dot{\delta\theta}} + \widehat{\omega} \widehat{\delta\theta}  - \widehat{\delta\theta} \widehat{\omega}
		\end{eqnarray}
		In vector form:
		\begin{eqnarray}
		\delta\omega = \dot{\delta\theta} + \widehat{\omega}\delta\theta =  \dot{\delta\theta} + \omega \times \delta\theta
		\end{eqnarray}
		Accordingly, we have the following:
		\begin{eqnarray}
			D_{\omega} \mathcal{L} \cdot \delta\omega &=& \omega^T J \delta\omega = \omega^T J \left(\dot{\delta\theta} + \widehat{\omega}\delta\theta\right) \\
			D_T\mathcal{L} \cdot\delta T &=& - \frac{\partial U}{\partial T}^T \delta T = - \frac{\partial U}{\partial T}^T T \widehat{\delta\theta} = - \frac{1}{2}\left[ \frac{\partial U}{\partial T}^T T - T^T \frac{\partial U}{\partial T} \right]_{\times}\cdot\delta\theta 
		\end{eqnarray}
		We can write the variation above as:
		\begin{equation}
		\delta S = \int_0^T \left[-\frac{1}{2}\left[ \frac{\partial U}{\partial T}^T T - T^T \frac{\partial U}{\partial T} \right]_{\times}\cdot\delta\theta + \omega^T J \left(\dot{\delta\theta} + \widehat{\omega}\delta\theta\right)\right] dt
		\end{equation}
		Rearranging terms we get:
		\begin{eqnarray}
		 \delta S = \int_0^T \left[\left(\widehat{\omega} J \omega - \frac{1}{2} \left[\frac{\partial U}{\partial T}^T T - T^T \frac{\partial U}{\partial T} \right]_{\times}\right)\cdot\delta\theta + \omega^T J \dot{\delta\theta} \right] dt
		\end{eqnarray}
		Performing integration by parts for the last term in the integrand we obtain:
		\begin{eqnarray}
		\delta S = \int_0^T \left(J \dot{\omega} + \widehat{\omega} J \omega - \frac{1}{2} \left[\frac{\partial U}{\partial T}^T T - T^T \frac{\partial U}{\partial T} \right]_{\times}\right)\cdot\delta\theta \ dt + \left.\omega^T J \delta\theta \right|_0^T = 0
		\end{eqnarray}
		Since we have fixed boundary conditions, the variations on the initial and final times vanish and we're left with the integration. For the remaining term to be equal to zero for all arbitrary $\delta\theta$ the integrand must be equal to zero, which is a result of the fundamental theorem of calculus of variation. \\
		Thus we have
		\begin{equation}
		J \dot{\omega} + \omega \times J \omega = \frac{1}{2} \left[\frac{\partial U}{\partial T}^T T - T^T \frac{\partial U}{\partial T} \right]_{\times}
		\end{equation}
		
		Which is the second equation of motion governing the rotational motion of a rigid body. Thus the system of equations governing the motion of rigid bodies evolving on \textbf{SO(3)} is:
		
		\begin{align}
		\dot{T} &= T \widehat{\omega} \\
		J \dot{\omega} + \omega \times J \omega &= \frac{1}{2} \left[\frac{\partial U}{\partial T}^T T - T^T \frac{\partial U}{\partial T} \right]_{\times}
		\end{align}
		
		In the case of the presence of a non-conservative force, we instead apply D'Alembert principle:
		\begin{eqnarray}
		 \delta S + \int_0^T \delta W \ dt = 0 
		\end{eqnarray}
		The virtual work is defined as:
		\begin{eqnarray}
		\delta W = M^T \delta\theta
		\end{eqnarray}
		where $M$ is the non-conservative moment. Using the same variation of the action integral as before and combining terms together we get:
		
		\begin{align}
		\dot{T} &= T \widehat{\omega} \\
		J \dot{\omega} + \omega \times J \omega &= M + \frac{1}{2} \left[\frac{\partial U}{\partial T}^T T - T^T \frac{\partial U}{\partial T} \right]_{\times}
		\end{align}
		\newpage
		\section{Discrete Variational Mechanics}
		\subsection{Derivation of the Symplectic Variational Integrator Scheme}
		In the previous sections, we formulated the properties of rotation group \textbf{SO(3)} and derived the equations of rotation motion of a rigid body in the continuous case. The continuous equations of motion are non-linear differential equations and have no closed form solutions except for trivial cases, and thus they are often solved numerically. However, a direct discretization of the differential equations using normal integrators, such as the 4th order Runge-Kutta, does not necessarily preserve the symplectic structure that the phase space is endowed with, and thus such integrators may lead to erroneous results. In this chapter we develop a 2nd order variational integrator for the rigid body motion from the variational principle of D'Alembert. 
		
		We begin by the same definition of the continuous Lagrangian:
		$$ \mathcal{L} = \frac{1}{2}\omega^T J \omega - \textbf{U}(T) $$
		and the action integral:
		$$ \textbf{S} = \int_0^T \mathcal{L} \ dt $$
		
		We define a discrete Lagrangian and use the mid point rule to approximate the action integral along the path of the system as follows:
		$$\mathcal{L}_{d,k+\frac{1}{2}}  = \frac{1}{2}\omega_{k+\frac{1}{2}}^T J \omega_{k+\frac{1}{2}} - \textbf{U}(T_{k+\frac{1}{2}}) $$
		$$ \textbf{S} \approx \textbf{S}_d = \sum_{k=0}^{N-1} \mathcal{L}_{d,k+\frac{1}{2}} \Delta$$
		where $\Delta$ is the time step of the discretization scheme, and $t = k\ \Delta$ is the current discrete time.
		
		We use the following midpoint approximations for the system states as follows:
		\begin{equation} \label{eq:tmean} T_{k+1} + T_{k} = V \ \kstep{T} \end{equation}
		where in Eq.(\ref{eq:tmean}) we utilize a result given in \cite{moakher}, which simply states that the mean of two rotation matrices is equivalent to the polar decomposition of their standard euclidean mean, and we have:
		$$ V = \sqrt{\left(T_{k+1} + T_k\right) \left(T_{k+1} + T_k\right)^T} $$
		And we define the matrix of relative transformation between two orientations at time steps $k+1$ and $k$ as $\kstep{R}$, and use Rodriguez formula for computing the corresponding rotation vector that represents this small rotation:
		\begin{equation} \label{eq:omega} \kstep{R}=  T_{k+1} T_{k}^T = \mathbb{I}_{3\times 3} + \frac{\sin(\kstep{\psi})}{\kstep{\psi}} \kstep{\widehat{\psi}} + \frac{1- \cos(\kstep{\psi})}{\kstep{\psi^2}} \kstep{\widehat{\psi}^2} \end{equation}
		And approximating the angular velocity vector as: 
		$$ \kstep{\omega} = \frac{1}{\Delta} \kstep{T}^T \kstep{\psi} $$
		Where we had to multiply by the inverse of the mean rotation $\kstep{T}$ because the way we defined the matrix $\kstep{R}$ will give us the rotation vector in space coordinates, where as we want the angular velocity in body coordinates. \\
		
		Applying the discrete version of D'Alembert principle:
		\begin{equation} \label{eq:dalembr} \delta \textbf{S} + \int_0^T \delta W\ dt \approx \delta \textbf{S}_d + \sum_{k=0}^{N-1} \kstep{\delta W} \Delta =  \sum_{k=0}^{N-1} \delta \mathcal{L}_d \Delta + \sum_{k=0}^{N-1} \kstep{\delta W} \Delta  = 0\end{equation}
		In order to proceed, we need to express the variations of the action integral and the virtual work in terms of $T_{k+1}$ and $T_k$.
		
		We know from the properties of \textbf{SO(3)} that
		\begin{eqnarray}
		\delta T = \widehat{\delta\theta} T
		\end{eqnarray}
		Thus computing the variation of Eq.\ref{eq:tmean}, we obtain
		\begin{eqnarray}
		\delta T_{k+1} + \delta T_{k} &=& \delta V \kstep{T} + V \delta \kstep{T} \\
		 \widehat{\delta\theta}_{k+1} T_{k+1} + \widehat{\delta\theta}_{k} T_{k} &=& \delta V \kstep{T} + V \kstep{\widehat{\delta\theta}}\kstep{T} \\
		 \widehat{\delta\theta}_{k+1} T_{k+1} \kstep{T}^T + \widehat{\delta\theta}_{k} T_{k} \kstep{T}^T &=& \delta V + V \kstep{\widehat{\delta\theta}}
		 \end{eqnarray}
		Define:
		\begin{eqnarray}
		Y_{k+1} &=& T_{k+1} \kstep{T}^T \\
		Y_k &=& T_{k} \kstep{T}^T
		\\
		 \delta T_{k+1} + \delta T_{k} &=& \delta V \kstep{T} + V \delta \kstep{T} \\
		 \widehat{\delta\theta}_{k+1} T_{k+1} + \widehat{\delta\theta}_{k} T_{k} &=& \delta V \kstep{T} + V \kstep{\widehat{\delta\theta}}\kstep{T} \\
		 \widehat{\delta\theta}_{k+1} T_{k+1} \kstep{T}^T + \widehat{\delta\theta}_{k} T_{k} \kstep{T}^T &=& \delta V + V \kstep{\widehat{\delta\theta}} 
		\end{eqnarray}
		Define:
		\begin{eqnarray}
		 Y_{k+1} &=& T_{k+1} \kstep{T}^T 
		 Y_k = T_{k} \kstep{T}^T \\
		 Y_{k+1} &=& T_{k+1} \kstep{T}^T \\
		 Y_k &=& T_{k} \kstep{T}^T 
		\end{eqnarray}
		We can write:
		\begin{eqnarray}
	\widehat{\delta\theta}_{k+1} Y_{k+1} + \widehat{\delta\theta}_{k} Y_k = \delta V + V \kstep{\widehat{\delta\theta}}
		\end{eqnarray}
		Taking anti-symmetric part followed by the hodge star operator, and noting out that $\delta V$ is a symmetric matrix, we end up with:
		\begin{eqnarray}
		\widetilde{Y}_{k+1}^T \delta\theta_{k+1} + \widetilde{Y}_{k}^T \delta\theta_{k} = \widetilde{V}\delta\kstep{\theta}
		\end{eqnarray} 
		where we have used a property of the hat map to obtain: 
		\begin{eqnarray}
		 \widetilde{A} = Tr\left(A\right)\cdot\mathbb{I}_{3\times3} - A 
		\end{eqnarray}
		Thus:
		\begin{equation} \label{eq:var_tm}\delta\kstep{\theta} =\widetilde{V}^{-1} \widetilde{Y}_{k+1}^T \delta\theta_{k+1}  + \widetilde{V}^{-1} \widetilde{Y}_{k}^T \delta\theta_{k} \end{equation}
		Taking the anti symmetric part of equation \ref{eq:omega}, we get:
		$$d \kstep{R} = \delta T_{k+1} T_{k} + T_{k+1} \delta T_{k}^\intercal$$
		$$\kstep{\widehat{\delta\phi}} \kstep{R} = \widehat{\delta\theta}_{k+1} T_{k+1} T_k^\intercal + T_{k+1} \left[\widehat{\delta\theta}_{k} T_{k}\right]^\intercal$$
		$$\kstep{\widehat{\delta\phi}} \kstep{R} = \widehat{\delta\theta}_{k+1} T_{k+1} T_k^\intercal - T_{k+1} T_{k}^\intercal \widehat{\delta\theta}_{k} $$
		$$\kstep{\widehat{\delta\phi}} \kstep{R} = \widehat{\delta\theta}_{k+1} \kstep{R} - \kstep{R} \widehat{\delta\theta}_{k}$$
		$$\kstep{\widehat{\delta\phi}} = \widehat{\delta\theta}_{k+1} - \kstep{R} \widehat{\delta\theta}_{k} \kstep{R}^\intercal $$
		$$\kstep{\delta\phi} = \delta\theta_{k+1} - \kstep{R} \delta\theta_{k}$$
		$$ \frac{1}{2} \left[\delta \kstep{R} - \left(\delta\kstep{R}\right)^\intercal\right]^\vee = \frac{\norm{\kstep{\psi}}\cos\norm{\kstep{\psi}} - \sin\norm{\kstep{\psi}}}{\norm{\kstep{\psi}}^2} \delta\norm{\kstep{\psi}} \kstep{\psi} + \frac{\sin\norm{\kstep{\psi}}}{\norm{\kstep{\psi}}} \delta\kstep{\psi}$$
		$$\norm{\kstep{\psi}}^2 = \kstep{\psi}^\intercal\kstep{\psi}$$
		$$\norm{\kstep{\psi}} \cdot \delta\norm{\kstep{\psi}} = \kstep{\psi}^\intercal \delta\kstep{\psi} $$
		$$\delta\norm{\kstep{\psi}} = \frac{1}{\norm{\kstep{\psi}}} \kstep{\psi}^\intercal \delta\kstep{\psi}$$
		$$ \frac{1}{2} \left[\delta  \kstep{R} - \left(\delta\kstep{R}\right)^\intercal\right]^\vee = \left[\frac{\norm{\kstep{\psi}}\cos\norm{\kstep{\psi}} - \sin\norm{\kstep{\psi}}}{\norm{\kstep{\psi}}^3} \kstep{\psi} \kstep{\psi} ^\intercal+ \frac{\sin\norm{\kstep{\psi}}}{\norm{\kstep{\psi}}}\mathbb{I}_{3\times3}\right] \delta\kstep{\psi}$$
		$$\left[\delta  \kstep{R} - \left(\delta\kstep{R}\right)^\intercal\right]^\vee = \left[\kstep{\widehat{\delta\phi}} \kstep{R} + \kstep{R}^\intercal\kstep{\widehat{\delta\phi}}\right]^\vee = \kstep{\widetilde{R}} \kstep{\delta\phi}  $$
		$$ \left[\frac{\norm{\kstep{\psi}}\cos\norm{\kstep{\psi}} - \sin\norm{\kstep{\psi}}}{\norm{\kstep{\psi}}^3} \kstep{\psi} \kstep{\psi} ^\intercal+ \frac{\sin\norm{\kstep{\psi}}}{\norm{\kstep{\psi}}}\mathbb{I}_{3\times3}\right] \delta\kstep{\psi} = \kstep{\widetilde{R}} \kstep{\delta\phi}$$
		$$\kstep{F}\delta\kstep{\psi} = \frac{1}{2} \kstep{\widetilde{R}}\delta\kstep{\phi}$$
		$$\delta\kstep{\psi} = \frac{1}{2}\kstep{F}^{-1}\kstep{\widetilde{R}}\delta\kstep{\phi}$$
		$$\delta\kstep{\psi} = \frac{1}{2}\kstep{F}^{-1}\kstep{\widetilde{R}}\delta\theta_{k+1} - \frac{1}{2}\kstep{F}^{-1}\kstep{\widetilde{R}} \kstep{R}\delta\theta_{k}$$
		$$\delta\kstep{\omega} \; \Delta = \left[\delta \kstep{T}\right]^\intercal \kstep{\psi} + \kstep{T}^\intercal \delta\kstep{\psi} $$
		$$\delta\kstep{\omega} \; \Delta = -\kstep{T}^\intercal \kstep{\widehat{\delta\theta}} \kstep{\psi}+ \kstep{T}^\intercal \delta\kstep{\psi} $$
		$$\delta\kstep{\omega} \; \Delta =  \kstep{T}^\intercal \kstep{\widehat{\psi}} \kstep{\delta\theta}+ \kstep{T}^\intercal \delta\kstep{\psi} $$
		$$\delta\kstep{\omega} \; \Delta =  \kstep{T}^\intercal \left[\kstep{\widehat{\psi}} \kstep{\delta\theta}+ \delta\kstep{\psi}\right] $$
		Expanding and rearranging:
		$$\delta\kstep{\omega} \; dt = \kstep{T}^\intercal\left[ \kstep{\widehat{\psi}}\widetilde{V}^{-1} \widetilde{Y}_{k} - \frac{1}{2}\kstep{F}^{-1}\kstep{\widetilde{R}} \kstep{R}\right]\delta\theta_{k} \; + \;  \kstep{T}^\intercal \left[\kstep{\widehat{\psi}}\widetilde{V}^{-1} \widetilde{Y}_{k+1} + \frac{1}{2}\kstep{F}^{-1}\kstep{\widetilde{R}}\right]\delta\theta_{k+1}$$
		
		$$\mathcal{L}_d = \frac{1}{2} \kstep{\omega}^\intercal J \kstep{\omega} - \textbf{U}(\kstep{T})$$
		$$ \delta \int_0^T\mathcal{L}(\omega,T)\;dt = \int_0^T\delta\mathcal{L}(\omega,T)\;dt \approx \sum_{k=0}^{N-1} \delta\mathcal{L}_d(\kstep{\omega},\kstep{T}) \; \Delta $$
		$$\delta \mathcal{L}_d(\kstep{\omega},\kstep{T}) = \kstep{\omega}^\intercal J \delta\kstep{\omega} + \delta\textbf{U}(\kstep{T}) $$
		$$ \delta\textbf{U}(\kstep{T}) = \delta\textbf{U}_1(T_k,T_{k+1}) \delta\theta_{k} + \delta\textbf{U}_2(T_k,T_{k+1})\delta\theta_{k+1} $$
		
		$$ \sum_{k=1}^{N} \delta \mathcal{L}_d(\kstep{\omega},\kstep{T}) \Delta  = \sum_{k=1}^{N}\left[\kstep{\omega}^\intercal J \delta\kstep{\omega} - \delta\textbf{U}_1(T_k,T_{k+1}) \delta\theta_{k} - \delta\textbf{U}_2(T_k,T_{k+1})\delta\theta_{k+1}\right] \Delta $$
		
		Substituting and rearranging we get:
		$$\sum_{k=1}^{N-1} \delta \mathcal{L}_d(\kstep{\omega},\kstep{T}) \Delta  = \Theta_0^+ \delta\theta_{0} - \Theta_N^-\delta\theta_{N} + \sum_{k=1}^{N-1} \left[D_1 \mathcal{L}_d(T_k,T_{k+1})+D_2\mathcal{L}_d(T_{k-1},T_k)\right] \cdot \delta\theta_{k} \Delta $$
		Where:
		\begin{align} 
		D_1 \mathcal{L}_d(T_k,T_{k+1}) &= \frac{1}{\Delta t} \kstep{\omega}^\intercal J \kstep{T}^\intercal\left[\kstep{\widehat{\psi}}\widetilde{V}^{-1} \widetilde{Y}_{k} - \frac{1}{2}\kstep{F}^{-1}\kstep{\widetilde{R}} \kstep{R}\right] -  \delta\textbf{U}_1(T_k,T_{k+1}) \\
		D_2 \mathcal{L}_d(T_{k-1},T_k) &=  \frac{1}{\Delta t} \kstem{\omega}^\intercal J \kstem{T}^\intercal\left[\kstem{\widehat{\psi}}\widetilde{V}^{-1} \widetilde{Y}_{k} + \frac{1}{2}\kstem{F}^{-1}\kstem{\widetilde{R}}\right] - \delta\textbf{U}_2(T_{k-1},T_k) 
		\end{align}
		Rearranging terms in the sum:
		\begin{eqnarray}
		\sum_{k=1}^{N} \delta \mathcal{L}_d(\kstep{\omega},\kstep{T}) \Delta  = \sum_{k=1}^{N-1} \left[\Theta_k^+ - \Theta_k^-\right]\delta\theta_{k} \Delta   + \Theta_{0}^+ \delta\theta_{0}\Delta - \Theta_{N}^- \delta\theta_{N}\Delta
		\end{eqnarray}
		As for the variation of the virtual work, we have:
		\begin{eqnarray}
		\int_0^T \delta W \ dt \approx \sum_{k=0}^{N-1} \kstep{\delta W} \Delta = \sum_{k=0}^{N-1} \kstep{M} \cdot \kstep{\delta\theta} \Delta  = \sum_{k=0}^{N-1} \kstep{M} \cdot \left[\widetilde{V}^{-1} \widetilde{Y}_{k+1}^T \delta\theta_{k+1}  + \widetilde{V}^{-1} \widetilde{Y}_{k}^T \delta\theta_{k}\right] \Delta
		\end{eqnarray}
		Where, $\kstep{M}$ is the non-conservative moment acting on the rigid body. \\
		Define:
		\begin{eqnarray}
 \mathbb{F}^+_k &=& \kstem{M} \cdot \kstem{\widetilde{V}}^{-1} \widetilde{Y}_{k} \\
 \mathbb{F}^- _k &=& \kstep{M} \cdot \kstep{\widetilde{V}}^{-1} \widetilde{Y}_{k} 
		\end{eqnarray}
		Define:
		\begin{eqnarray}
		 \Theta_{k}^- = -D_1 \mathcal{L}_d(T_k,T_{k+1})\\
		 \Theta_{k}^+ = D_2 \mathcal{L}_d(T_{k-1},T_{k})
		\end{eqnarray}
		Substituting all of the above into Eq\ref{eq:dalembr} and after some rearranging of the terms, we get:
		\begin{equation}
		\sum_{k=1}^{N-1} \left[\Theta_k^+ - \Theta_k^- + F_k^+ + F_k^-\right] \delta\theta_k \Delta 
		\end{equation}
		For the above variation to be equal to zero, for all $\delta\theta_k$ which are arbitrary and independent, the expression in brackets mush vanish, which gives the discrete version of Lagrange equations of motion:
		\begin{equation}
		\label{eq:disc}
		\Theta_k^+ - \Theta_k^- + F_k^+ + F_k^- = 0
		\end{equation}
		Eq.\ref{eq:disc} is the rule that defines the propagation rule for the discrete dynamics of the rotational degrees of freedom of a rigid body. Given $(T_0,\omega_0)$ we solve the above equation for $(T_1,\omega_1)$ and thus the repetitive solution of this equation defines a map $\mathbb{\phi}_d $
		$$
		\mathbb{\phi}_d : (T_k,\omega_k) \rightarrow (T_{k+1},\omega_{k+1})
		$$
		which is the discrete flow of the discrete forced Lagrangian vector field. 
		
		\newpage
		\section{Local Stability of Nonlinear Systems}
			The concept of stability as we know it for dynamical systems was explored with the theory of stability by Lyapunov, a Russian mathematician and physicist, in \textit{The General Problem of the Stability of Motion}, 1892 \cite{lyapunov}. \\
			Consider the following autonomous system, no input included, 
			\begin{eqnarray}
			\dot{x} =f(x)
			\end{eqnarray}
			with $x_0$ is its equilibrium point, $f(x_0)=0$. \\
			$x$ represents the local coordinates on the configuration manifold, $f(x)$ is the drift smooth vector field on the configuration manifold. \\
			
			The general concept of stability declares that $x_0$ is a locally stable equilibrium point if for any neighborhood $\epsilon$ of $x_0$, there exists a neighborhood $\delta$ such that if the set of initial conditions $\tilde{x} \in \delta$ , the orbits/solutions $x(t,0,\tilde{x}) \in \epsilon$ for all $t$. \\
			$x_0$ is considered a locally asymptotically stable equilibrium, if it is locally stable and there exists a neighborhood $\delta_0$ such that if $\tilde{x} \in \delta_0$, the orbits/solutions $x(t,0,\tilde{x}) \longrightarrow x_0$ as $t\longrightarrow \infty$.
			\subsection{First Method of Lyapunov}
			Consider the following linearized system,
			\begin{eqnarray}
			\dot{x} = A \ x
			\end{eqnarray}
			where $A = \frac{\partial f}{\partial x} |_{x_0}$ \\
			The first method of Lyapunov states that $x_0$ is considered a locally asymptotically stable equilibrium if all eigenvalues of $A \in C^-$. $x_0$ is considered a locally unstable equilibrium if at at least one of the eigenvalues of $A \in C^+$. \\
			The method is not decisive about the stability of the equilibrium point if one of the eigenvalues is located on the imaginary axis. 
			\subsection{Second Method of Lyapunov and LaSalle-Yoshizawa Theorem}
			The second or the direct method of Lyapunov aims at deciding on  the stability of an equilibrium point from checking some properties of a defined smooth positive definite function $V$, called Lyapunov function, and its Lie derivative along the system dynamics. \\
			Consider the smooth Lyapunov function $V$ defined on the neighborhood $\delta_0$ where
			\begin{eqnarray}
			V(x_0) &=&0  \\
			V(x) &>&0, \ \ \ \ \ x \neq x_0
			\end{eqnarray}
			Before we relate the properties of the Lyapunov function with the system, we need to define an invariance of a set:
			A set $E$ in a manifold $M$ is invariant if for all the initial conditions $\tilde{x} \in E$, all the solutions/orbits $x(t,0,\tilde{x})$ stay in $E$ for all $t$.\\
			
			The second method of Lyapunov dictates $x_0$ is a locally stable equilibrium point if $L_f V(x) \leq 0$ for all $\tilde{x} \in \delta_0$. \\
			
			Let $E_0 = \{x \in \delta_0 \ | \ L_f V(x) =0 \}$ and $E$ is the largest invariant set $\in E_0$ where
			every solution/orbit $x(t,0,\tilde{x})$ starting in $\tilde{x} \in E$ converges to $E$. \\
			To guarantee asymptotic stability, we would need to strictly condition that every solution/orbit $x(t,0,\tilde{x})$ starting in $E$ converges to $x_0$ i.e. we need to make sure that the largest invariant set under the dynamics is the equilibrium point itself. In other words, we only permit the equilibrium point to be have $L_f V(x) =0$. This condition is known as LaSalle-Yoshizawa Invariance Theorem.
			\subsection{Stabilization of A Nonlinear Control  System}
			Consider the following control affine system
			\begin{eqnarray}
			\dot{x} =f(x)+ \sum_{i}^{m} g_i(x) u_i
			\end{eqnarray}
			with $x_0$ is its equilibrium point, $f(x_0,u_0)=0$. \\
			$x$ represents the local coordinates on the state space manifold, $f(x), g(x)$ are the respective drift and control smooth vector fields on the state space manifold, $u \in \mathbb{R}^m$.  \\
			
			Assume a smooth feedback control $u = \eta(x)$. We choose $\eta(x)$ such that for $V(x)$ defined in $\delta_0$
			\begin{eqnarray}
			\eta_i (x) = -L_{g_{i}} V(x)
			\end{eqnarray}
			The condition of stability is consequently,
			\begin{eqnarray}
			L_f V(x)+\sum_{i}^{m} L_{\eta_i g_i} V(x) \leq 0 \\
			L_f V(x)-\sum_{i}^{m} (L_{g_i} V(x))^2 \leq 0 
			\end{eqnarray}
			let $E_0 = \{x \in \delta_0 \ | \ L_f V(x)-\sum_{i}^{m} (L_{g_i} V(x))^2 =0 \}$. However, we need to make sure that the solutions emerged from uncontrolled dynamics over time converge to the equilibrium point. We define $E_0$ as $E_0 = \{x \in \delta_0 \ | \ L_f V(x) =0, \ \sum_{i}^{m} (L_{g_i} V(x))^2 =0 \}$ and $E$ is the largest invariant set under the dynamics $\in E_0$. \\
			To guarantee asymptotic stability, we would need to strictly condition that every solution/orbit $x(t,0,\tilde{x})$ starting in $E$ and converges to $x_0$ i.e. we need to ensure that the largest invariant set under the drift dynamics is the trivial solution which is the equilibrium point. Owing to the fact that any solution in $E_0$ dictates $\sum_{i}^{m} (L_{g_i} V(x))^2 =0$. \\
			
			If $x_0$ is asymptotically stable, the set of initial conditions  $\tilde{x} \in \delta$ satisfying the orbits/solutions $x(t,0,\tilde{x})\longrightarrow x_0$ as $t\longrightarrow \infty$ is called the region of attraction. If the region of attraction is the whole manifold, $x_0$ is a global asymptotically stable equilibrium.
			\section{Backstepping Control}
			Developed first by Kokotovic in \textit{The joy of feedback: nonlinear and adaptive}, 1992 \cite{Kokotovic}.
			Backstepping is a type of controller made especially to exploit the recursive structures of some nonlinear dynamical systems such as
			\begin{equation}
			\begin{aligned}
			\dot{x}_1 &= f(x_1) +g(x_1) x_2 \\
			\dot{x}_2 &= f(x_1,x_2) +g(x_1,x_2) x_3 \\
			\vdots \\
			\dot{x}_k &=f(x_1,x_2,..,x_k) +g(x_1,x_2,..,x_k) u
			\end{aligned}
			\end{equation} 
			$x$ represents the local coordinates on the state space manifold, $u \in \mathbb{R}$.  \\
			Notice that the control input only appears in one state equation. Therefore, the controller is designed via stabilizing the first state assuming $x_2$ is a virtual control input. After that, $x_3$ progressively stabilizes
			the second state equation to follow the virtual control of the first. Recursively, $x_k$ virtually stabilizes $x_{k-1}$. The objective is to acquire the actual control input $u$ to stabilize the whole system. This is achieved with the help of Lyapunov’s Direct method of stability.\\
			Consider the following single integrator system, 
			\begin{equation}
			\begin{aligned}
			\dot{x}_1 &= f(x_1) +g(x_1) x_2 \\
			\dot{x}_2 &= f(x_1,x_2) +g(x_1,x_2) u \\
			\end{aligned}
			\end{equation}
			We would choose a Control Lyapunov Function (clf) as follows to fulfill asymptotic stability for the first state:
			\begin{eqnarray}
			V(x_1) =\frac{1}{2} x_1^2
			\end{eqnarray}
			This requires that the lie derivative of clf along the dynamics must be negative $L_{f(x_1)}V(x_1) + L_{g(x_1)}V(x_1) <0$. This step derives an expression of the targeted $x_2$ to stabilize the first equation. Augmenting, the clf to
			be a function of the error between targeted $x_2$ and its current value, the control input will be proportional to this error bringing about the asymptotic stability of the whole system. 
			\begin{eqnarray}
			V_a(x_1,x_2) =\frac{1}{2} x_1^2 + \frac{1}{2} (x_2 -x_{2,tar})^2 
			\end{eqnarray}
			where $x_{2,tar}$ is the fictitious control of the first state equation. \\ 
			Once again, the lie derivative of the augmented clf along the dynamics must be negative which yields an expression of the actual control u.
			\subsection{Application of Backstepping Control In Rotational Motion In SO(3)}
			The most profound privilege of geometric nonlinear control is the capability of preserving global nature of the attitude dynamics. Therefore, The stabilization or tracking is enabled even in cases of large attitude errors like almost $180^ \circ$ rotations. \\
			Consider the following rigid body rotation EOM:
			\begin{align}
			\dot{R} &= R \ \widehat{\Omega} \\
			J \ \dot{\Omega} &= - \Omega \times J \Omega + \ q 
			\end{align}
			With desired rotation matrix $R_d$ and angular velocity $\Omega_d$ must be fulfilled using certain control input $q$. where\\
			\begin{eqnarray}
			\dot{R}_d &=& R_d \ \widehat{\Omega}_d
			\end{eqnarray}
			Since the configuration space of the rigid body attitude  is $SO(3)$, it is essential to choose an error function consistent with the structure of the rotation group. Hence, we define the Riemannian distance/metric between the desired and the actual  rotation through multiplication, the group operation defined on SO(3) manifold,
			\begin{eqnarray}
			E = R_d^T R   
			\end{eqnarray}
			where $E \in \textbf{SO(3)}$ \\
			The error function in rotation is handled by defining a function of the trace of $E$.
			\begin{eqnarray}
			\psi(R, R_d) = \psi(tr(E))
			\end{eqnarray}
			where $\psi: \textbf{SO(3)} \times \textbf{SO(3)} \rightarrow \mathbb{R}$. \\
			Here, the error function is chosen the same as \cite{Lee}:
			\begin{eqnarray} 
			\psi(R, R_d) = 2- \sqrt{1+tr(E)}
			\end{eqnarray} 
			Notice the error function in rotation is positive definite since $-1 \leq tr(A)\leq 3$ for any matrix $A$. Moreover, the trace of any rotation matrix can be said to equal $1+2\ \cos(\theta)$ that is why the rotation error is a quadratic function in $\theta$. Consequently, a clf candidate for rotation can be the error function itself.\\
			
			A rotation error vector $e_R$ is defined for the error function $\psi$ at any point $p$ as $e_R \in T_p \textbf{SO(3)}$. where
			\begin{eqnarray}
			e_R = \frac{1}{2 \sqrt{1+tr(E)}}\ (E-E^T)^\vee
			\end{eqnarray}
			where $(E-E^T)^\vee$ is the Hodge dual to the skew-symmetric matrix $(E-E^T)$.\\
			\textbf{Proof}: \\
			Following the upcoming two statements:
			\begin{eqnarray}
			\delta \psi(R,R_d) &=& - \frac{1}{2 \sqrt{1+tr(R_d^T R)}} \ \delta[tr(R_d^T R)] \\
			\delta[tr(R_d^T R)] &=& tr(\delta[R_d^T R]) \\
			\end{eqnarray}
			The variation in the Riemannian metric is expressed as follows:
			\begin{eqnarray}
			\delta E = R_d^T\ \ \delta R \\
			\delta E = R_d^T\ R\ \widehat{\partial R} 
			\end{eqnarray}
			where the element of Lie algebra is $\widehat{\partial R}$.\\
			Accordingly,
			\begin{eqnarray}
			\delta E &=& E\ \widehat{\partial R} \\
			\delta[tr(R_d^T R)] &=& tr( E\ \widehat{\partial R}) \\
			\end{eqnarray}
			Equivalently,
			\begin{eqnarray}
			\delta [tr(R_d^T R)] &=& - (R_d^T R - R^T R_d)^\vee \cdot \ {\partial R} \\
			\delta [tr(E)] &=& (E - E^T)^\vee \cdot \ \widehat{\partial R}
			\end{eqnarray}
			Hence,
			\begin{eqnarray}
			\delta \psi(R,R_d) &=& \frac{\partial \psi(R,R_d)}{\partial tr(R_d^T R)} \ 	\delta [tr(R_d^T R)] \cdot \ \widehat{\partial R} \\
			e_R &=& \frac{1}{2 \sqrt{1+tr(E)}}\ (E-E^T)^\vee
			\end{eqnarray}
			From Rodriguez' Formula,
			\begin{eqnarray}
			E &=& I+ \ \frac{\sin \theta}{\theta} \hat{\theta}+ \frac{1-\cos\theta}{\theta^2} \hat{\theta}^2 \\
			(E-E^T)^\vee &=& \frac{\sin \theta}{\theta} \hat{\theta} \\
			\psi(R,R_d) &=& 4 \ {\sin^2 \frac{||\theta||}{4}} \\
			||e_R||^2 &=& \sin^2 \frac{||\theta||}{2}
			\end{eqnarray}	
			Therefore, the rotation error vector is  in the direction of the Euler axis. The magnitude of the error vector has zero value if the rotation about the Euler axis is $0^\circ$ . 
			\begin{figure}[H]
				\centering
				\includegraphics[scale =.8]{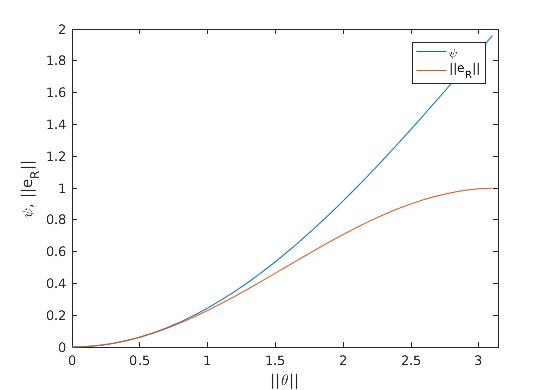}
				\caption{change of $\psi, ||e_R||$ with $||\theta||$.}
				\centering
			\end{figure}
			Basically, to apply the backstepping technique, we define
			\begin{eqnarray}
			V_R = \psi(R,R_d)
			\end{eqnarray}
			The error in the angular velocity represented in the body fixed frame is 
			\begin{eqnarray}
			e_\Omega = \Omega -R^T R_d \Omega_d
			\end{eqnarray}
			The Lie derivative of the clf along the error dynamics must be negative to reach asymptotic stability.
			\begin{eqnarray}
			\dot{V}_R &=& \frac{d }{dt} \psi(R,R_d) \\
			\frac{d }{dt} \psi(R,R_d) &=& - \frac{1}{2 \sqrt{1+tr(R_d^T R)}} tr[R_d^T R \widehat{\Omega} - \widehat{\Omega}_d R_d^T R ] \\
			R\ R_d^T \widehat{\Omega}_d R_d^T R &=& (R\ R_d^T \Omega_d)^\wedge \\
			\frac{d }{dt} \psi(R,R_d) &=& - \frac{1}{2 \sqrt{1+tr(R_d^T R)}} tr[R_d^T R \widehat{e}_\Omega]
			\end{eqnarray}
			Using the identity,
			\begin{eqnarray}
			tr[R_d^T R \widehat{e}_\Omega] = -(R_d^T R - R^T R_d)^\vee \cdot \ {e}_\Omega
			\end{eqnarray}
			The Lie derivative of the clf of rotation:
			\begin{eqnarray}
			\dot{V}_R &=& \frac{1}{2 \sqrt{1+tr(R_d^T R)}} \ (R_d^T R - R^T R_d)^\vee \cdot \ {e}_\Omega 
			\end{eqnarray}
			Equivalently,
			\begin{eqnarray}
			\dot{V}_R &=& e_R \cdot \ {e}_\Omega 
			\end{eqnarray}
			The virtual control to asymptotically stabilize the rotation equation is
			\begin{eqnarray}
			\Omega_{tar} &=& -P \ e_R + R^T R_d \ \Omega_d
			\end{eqnarray}
			where $P$ is any positive definite matrix. \\
			The augmented clf for the control system is therefore,
			\begin{eqnarray}
			V_a &=& k_R \ \psi+\frac{1}{2} (\Omega -\Omega_{tar}) \cdot S \ (\Omega -\Omega_{tar}) 
			\end{eqnarray}
			where $k_R$ is a positive constant and $S$ is any positive definite matrix. \\
			To ensure the system follows the targeted angular velocity, We assign $\dot{V}_a <0$.
			\begin{eqnarray}
			\dot{V}_a &=& k_R \frac{d}{dt} \psi +(\Omega -\Omega_{tar}) \cdot S \ (\dot{\Omega} -\dot{\Omega}_{tar}) \\
			\dot{\Omega}_{tar} &=& - \widehat{\Omega} \ R^T R_d \ \Omega_d +R^T R_d \ \widehat{\Omega}_d \ \Omega_d +R^T \ R_d \ \dot{\Omega}_d -P \ \dot{e}_R
			\end{eqnarray}
			Using the identity,
			\begin{eqnarray}
			\widehat{\Omega}_d \ \Omega_d = \Omega_d \times \Omega_d =0
			\end{eqnarray}
			Therefore,
			\begin{eqnarray}
			\dot{\Omega}_{tar} &=& - \widehat{\Omega} \ R^T R_d \ \Omega_d +R^T \ R_d \ \dot{\Omega}_d -P \ \dot{e}_R
			\end{eqnarray}
			$\dot{e}_R$ is obtained as follows:
			\begin{equation}
			\begin{split}
			\dot{e}_R = \frac{1}{2 \sqrt{1+tr(R_d^T R)}} \ (- \widehat{\Omega}_d\ R_d^T R +R_d^T R \ \widehat{\Omega}+ \widehat{\Omega} \ R^T R_d - R^T R_d \ \widehat{\Omega}_d)^\vee  \\ + (R_d^T R - R^T R_d)^\vee* \frac{-tr(R_d^T R \ \widehat{\Omega}-\widehat{\Omega}_d\ R_d^T R)}{4 (1+tr(R_d^T R))^\frac{3}{2}}
			\end{split}
			\end{equation}
			\begin{eqnarray}
			\dot{e}_R = \frac{1}{2 \sqrt{1+tr(R_d^T R)}} \ (R_d^T R \ \widehat{e}_\Omega+\widehat{e}_\Omega \ R^T R_d)^\vee - \frac{tr[R_d^T R \ \widehat{e}_\Omega]}{2 (1+tr(R_d^T R))} \ e_R
			\end{eqnarray}
			Using the following two identities,
			\begin{eqnarray}
			tr[R_d^T R \ \widehat{e}_\Omega] &=& - e_\Omega^T \ (R_d^T R - R^T R_d) \\
			(R_d^T R \ \widehat{e}_\Omega+\widehat{e}_\Omega \ R^T R_d)^\vee &=& (tr(R^T R_d) \mathbb{I} -R^T R_d) e_\Omega
			\end{eqnarray}
			We can derive the derivative of the error vector as:
			\begin{eqnarray}
			\dot{e}_R &=& \frac{1}{2 \sqrt{1+tr(R_d^T R)}}(2 e_R \ e_R^T +tr(R^T R_d) \mathbb{I} -R^T R_d) e_\Omega \\
			\dot{e}_R =\beta \ e_\Omega \  &,&  \ \beta =  \frac{1}{2 \sqrt{1+tr(R_d^T R)}}(2 e_R \ e_R^T +tr(R^T R_d) \mathbb{I} -R^T R_d)
			\end{eqnarray}
			Eventually,
			\begin{eqnarray}
			\dot{V}_a = k_R \frac{d}{dt} \psi +(\Omega -\Omega_{tar}) \cdot S \ (J^{-1} \ (-\Omega \times J \Omega +q) -\dot{\Omega}_{tar})
			\end{eqnarray}
			The control input needed for attitude tracking is found as:
			\begin{eqnarray}
			q &=& \Omega \times J \Omega + J \dot{\Omega}_{tar} -F (\Omega -\Omega_{tar}) \\
			q &=&  \Omega \times J \Omega + J R^T R_d \dot{\Omega}_d -J \ \widehat{\Omega} R^T R_d \ \Omega_d -J \beta e_\Omega -F \ (\Omega -\Omega_{tar})
			\end{eqnarray}
			where F can be any positive definite matrix. A convenient choice of $F, P$ might be the inertia matrix to give the weighting of the error components same as its corresponding moment of inertia.  \\
			As stated in \cite{Lee}, the controller proposed from the same error function achieves exponential stability with an attitude error close to $180^\circ$. 
			\subsection{Application of Backstepping Control In Translation Motion}
			Consider the translational motion of a Quad-rotor represented in the inertial frame of reference. 
			\begin{align}
			\dot{r} &= v \\
			m \ \dot{v} &= m G +R \ u_b
			\end{align}
			where $m$ is mass, $r,v$ are the respective position and velocity vectors, $r,v \in \mathbb{R}^3$, $G$ is the gravity vector $G = \begin{bmatrix}
			0 \\ 0 \\ -g
			\end{bmatrix}$ , $u_b$ is the thrust force in the body fixed frame.	\\
			The control input is required to fulfill a desired position $r_d$ in the inertial frame with a command $u_b$. Assume the error in position is 
			\begin{eqnarray}
			e_r &=& r- r_d \\
			V_r &=& \frac{1}{2} e_r \cdot A \ e_r \\
			\dot{V}_r &=& e_r \cdot A \ \dot{e}_r \\
			\dot{V}_r &=& (r- r_d) \cdot A\ (\dot{r} -\dot{r_d})
			\end{eqnarray}
			To achieve the asymptotic stability, we assign the following $\dot{r}_{tar}$ enabling $\dot{V}_r <0$.
			\begin{eqnarray}
			\dot{r}_{tar}=v_{tar}= \dot{r}_d -B \ e_r
			\end{eqnarray}
			where $B$ is any positive definite matrix. \\
			In order to acquire the control input for the tracking problem, we define the augmented Lyapunov function as:
			\begin{eqnarray}
			V_a = \frac{1}{2} e_r \cdot A \ e_r+ \frac{1}{2} (v-v_{tar}) \cdot C \ (v-v_{tar})
			\end{eqnarray}
			And its Lie derivative along the error dynamics is,
			\begin{eqnarray}
			\dot{V}_a &=& e_r \cdot A B \ e_r+(v-v_{tar}) \cdot C \ (\dot{v}-\dot{v}_{tar}) \\
			\dot{V}_a &=& e_r \cdot A B \ e_r+(v-v_{tar}) \cdot C \ (G+ \frac{1}{m} \ R \ u_b-\dot{v}_{tar})
			\end{eqnarray}
			where $C$ is any positive definite matrix. \\
			Therefore,
			\begin{eqnarray}
			R u_b = m \ \dot{v}_{tar} -m \ G -D(v-v_{tar})
			\end{eqnarray}
			where $D$ is any positive definite matrix.
			\subsection{Backstepping for a Quad-rotor}
			At this point, we can apply the geometric backstepping technique in any rigid body. This may find many applications such as satellite attitude control and full tracking of Quad-rotors. We choose the quad-rotor to apply backstepping technique for its simple configuration. It is an under-actuated system where there are only four control inputs:thrust and roll, pitch and yaw moments. As a consequence the four control inputs allow to track four outputs. \\
			We use the Quad-rotor configuration same as \cite{madani}.
				\begin{figure}[H]
					\centering
					\includegraphics[scale =1]{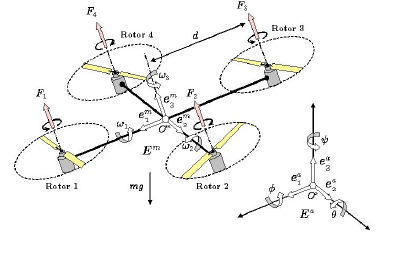}
					\caption{Quad-rotor Configuration \cite{madani}.}
					\centering
				\end{figure}
			In case of a quad-rotor, $f$ is the thrust force and it is in the positive direction of $e_3^m$. Hence, $u_b = f R e_3 $.\\
		   \begin{eqnarray}
		    f \ R \hat{e}_3 &=& m \ \dot{v}_{tar} -m \ G -D(v-v_{tar}) \\
		   f &=& [m \ \dot{v}_{tar} -m \ G -D(v-v_{tar})] \cdot R \hat{e}_3
		   \end{eqnarray}
		   We construct a matrix $R_c \in \textbf{SO(3)}$ where 
		   \begin{eqnarray}
		   R_c = [b_{1c}; b_{1c} \times b_{3c}; b_{3c}] \\
		   \end{eqnarray} 
		   and $b_{3c} = \frac{m \ \dot{v}_{tar} -m \ G -D(v-v_{tar})}{||m \ \dot{v}_{tar} -m \ G -D(v-v_{tar})||}$ , $b_{1c}$ is orthogonal to $b_{3c}$, $b_{1c},b_{3c} \in S^2$.\\
		   
		   $R_c$ represents the rotation required for the quad-rotor to follow a certain position command $r_d$.  $b_{3c}$ is taken in this direction to make sure the $z$ axis in the body frame represents the thrust vector for all time.
		   $R_c$ is constructed and then used as the desired rotation for the Quad-rotor for the attitude tracking in which it converges to this attitude with time. Therefore, we need to efficiently choose $b_{1c}$ for proper tracking performance. The user is thus required to enter a certain direction to fully define the matrix $b_{1c}$ and in turn $R_c$. \\
		   The closed loop system is presented as follows:
		   \begin{figure}[H]
		   	\centering
		   	\includegraphics[scale =.8]{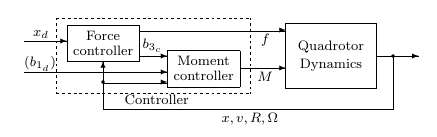}
		   	\caption{Quad-rotor Full Tracking \cite{leegeometric}.}
		   	\centering
		   \end{figure}   
		   For any input $b_{1d}$ not parallel to $b_{3c}$, $b_{1c}$ is defined as the projection of $b_{1d}$ on the plane perpendicular to $b_{3c}$:
		  \begin{eqnarray}
		  b_{1c} = -\frac{(b_{3c} \times (b_{3c} \times b_{1c}))}{||(b_{3c} \times (b_{3c} \times b_{1c}))||}
		  \end{eqnarray}
		  This method is proposed by Lee et al. in \cite{leegeometric}. We and apply the same technique with backstepping control (Results section). 
			\newpage
			\section{Aerodynamic Forces and Torques}
			Quad-rotor aerodynamic model is based on a combination of the two main theories: momentum theory and blade element theory. The two methods are used below to derive the aerodynamic forces and moments based on the work of Kroo et al. \cite{kroo}.
			\subsubsection{Momentum theory}
			Momentum theory is based on dealing with the rotor as an actuator disk  across which the flow is accelerated generating an inflow/induced velocity. Using conservation of mass and energy through disk, an expression for thrust can be derived $ T= 2 \rho A \nu_1 \sqrt{V^2 +\nu_1^2}$. \\ Solving for the inflow velocity leads to:
			\begin{eqnarray}
			\nu_1   = \left(\frac{V^2}{2}+ \sqrt{ \left(\frac{V^2}{2} \right)^2+\left(\frac{W}{2 \ \rho \ A}\right)^2 }\right)^{\frac{1}{2}}
			\end{eqnarray}
			where $\rho$ is the density of air, $W$ is rotor weight, $A$ is rotor disk area, $V$ is the horizontal velocity and $\nu$ is the inflow velocity. \\
			Momentum Theory assumes \cite{bangura}:
			\begin{description}
				\item [$\bullet$]infinite number of rotor blades hence a uniform constant force distribution is applied to rotor disc.
				\item [$\bullet$] very thin disc hence no resistance for air flow.
				\item [$\bullet$]irrational flow, no swirls.
				\item [$\bullet$] air outside control volume is undisturbed by the rotor disc.
			\end{description}
			Two important dimensionless quantities are always used in rotary literature: Inflow ratio and rotor advanced ratio.
			The inflow ratio relates the inflow velocity to the rotor tip velocity as follows:
			\begin{eqnarray}
			\lambda   =  \frac{\nu_1-\dot{z}}{\Omega R}
			\end{eqnarray}
			The rotor advanced ratio relates the sideways velocity to the rotor tip velocity as follows:
			\begin{eqnarray}
			\mu   =  \frac{V}{\Omega R}
			\end{eqnarray}
			Note that the sideways (horizontal) velocity is descried as \(V= \sqrt{\dot{x}^2+\dot{y}^2}\) and $\Omega$ is the angular velocity.
			\subsubsection{Blade Element Theory}
			Blade element theory is the method of determining the total aerodynamic forces and torques on a rotor by integrating the forces acting on single blade element 'airfoil' over the whole rotor. A demonstration of a blade element and local velocities and forces acting on it, is shown in next figure.
			\begin{figure}[H]
				\centering
				\includegraphics[scale =.3]{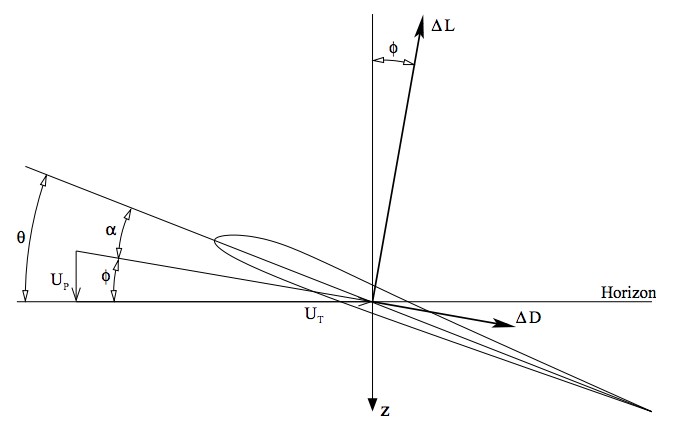}
				\caption{Blade Element\cite{kroo}}
				\centering
			\end{figure}
			Local velocity that is seen by the rotor is composed of two components: horizontal component due to the angular velocity of the element and its radial position and horizontal motion of the blade $U_T=\Omega R (\frac{r}{R}+ \mu \sin{\Psi})$,  and vertical component owing to the inflow and vertical motion of blade $U_P=\Omega R \lambda$. 
			Note that $\theta$ is the incidence angle, $\alpha$ is the angle of attack, $\Psi$ is the azimuth angle and $\phi$ is the inflow angle.
			
			\subsubsection{Thrust Force}
			As blade element theory proposes, thrust force can be obtained by integrating vertical forces applied to the airfoil over the whole rotor. The vertical forces acting on the airfoil are lift force component and drag force component defined as following: $ \Delta F_v = \Delta L \cos{\phi}-\Delta D \sin{\phi}$.\\
			Taken assumptions \cite{kroo}:
			\begin{description}
				\item [$\bullet$] Rotor blade has constant cord and entire rotor lies in one plane.
				\item [$\bullet$] Negligible aerodynamic moments i.e. sheer center and aerodynamic center are very close and stiff rotor.
				\item [$\bullet$] Negligible gravity torques due to rotor's light weight.
				\item [$\bullet$] Rigid blades therefore no blade flapping and coning.
				\item [$\bullet$]Coefficient of lift is linear in the angle of attack, $C_l=a \alpha= a(\theta-\phi)$.
				\item [$\bullet$] Linear twist distribution is used, $\theta=\theta_0 - \theta_{tw} (\frac{r}{R})$.
			\end{description}
			Lift and drag forces as functions of dynamic pressure $q= \frac{1}{2} \rho \ U_T^2$, reference area $S=\Bar{c} \ \Delta r $, $C_l$ and $\Bar{C_d}$.
			\begin{align}
			\Delta L   &=  \frac{1}{2} \rho \ U_T^2 \ a \Big(\theta_0 -\theta_{tw} \frac{r}{R}- \frac{U_P}{U_T}\Big)\\
			\Delta D   &=  \frac{1}{2} \rho \ U_T^2 \ \Bar{C_d} \ \Bar{c} \ \Delta r
			\end{align}
			Applying small angle approximation, vertical forces become \( \Delta F_v = \Delta L \) and integrating it to get the thrust force:
			\begin{align*}
			T   &= \frac{N}{2 \pi} \int_{0}^{2 \pi}\int_{0}^{R} \frac{\Delta L}{\Delta r} dr \ d\Psi\\
			&= N \rho \ a \ \Bar{c} \ (\Omega R)^2 R \ \Big[\Big(\frac{1}{6}+\frac{1}{4} \mu^2\Big) \theta_0 - (1+\mu^2) \frac{\theta_{tw}}{8}-\frac{\lambda}{4}\Big]  
			\end{align*}
			The coefficient of thrust is a dimensionless quantity defined by
			\begin{eqnarray}
			C_T= \frac{T}{\rho A (\Omega R)^2}
			\end{eqnarray} 
		Hence, 
			\begin{eqnarray}
			\frac{C_T}{\sigma a}   = \Big(\frac{1}{6}+\frac{1}{4} \mu^2\Big) \theta_0 - (1+\mu^2) \frac{\theta_{tw}}{8}-\frac{\lambda}{4}
			\end{eqnarray}
			This applies to all the aerodynamic coefficients derived later.
			\subsubsection{Hub Force}
			Similarly, to obtain the hub force, horizontal forces acting on the blade element must be integrated. The hub force have two components in x-direction where azimuth angle $\Psi=0$ called the H-force and the other force in the y-direction where $\Psi=\frac{\pi}{2}$ is the Y-Force. The  horizontal forces acting on blade element are as defined before. Components of lift and drag $\Delta F_H =\Delta L \sin{\phi}+\Delta D \cos{\phi}$, a small angle approximation changes the expression to $\Delta F_H =\Delta L +\Delta D \frac{U_P}{U_T}$.\\
			The H-force is: 
			\begin{align}
			H   &= \frac{N}{2 \pi} \int_{0}^{2 \pi}\int_{0}^{R} \Big[\frac{\Delta D}{\Delta r}+\frac{\Delta L}{\Delta r} \frac{U_P}{U_T}\Big]\sin{\Psi}  dr \ d\Psi\\
			&= N \rho a \ \Bar{c} \ (\Omega R)^2 R \Big[\frac{\mu \ \Bar{C_d}}{4a}+\frac{1}{4} \lambda \mu \Big(\theta_0 -\frac{\theta_{tw}}{2}\Big)\Big]\\
			\end{align}
			And the Coefficient of hub force:
			\begin{eqnarray}
			\frac{C_H}{\sigma a}   = \frac{\mu \ \Bar{C_d}}{4a}+\frac{1}{4} \lambda \mu \Big(\theta_0 -\frac{\theta_{tw}}{2}\Big)
			\end{eqnarray}
			Similarly in $\Psi=\pi/2$, the Y-force is: 
			\begin{eqnarray}
			Y   = -\frac{N}{2 \pi} \int_{0}^{2 \pi}\int_{0}^{R} \Big[\frac{\Delta D}{\Delta r}+\frac{\Delta L}{\Delta r} \frac{U_P}{U_T}\Big]\cos{\Psi} \  dr \ d\Psi 
			\end{eqnarray}
			The integration for Y-force equals to zero, and $C_Y$ is therefore also zero.
			\subsubsection{Torques}
			Similar to forces, aerodynamic torques are determined by integrating forces acting on blade element multiplied by moment arm over entire rotor. The aerodynamic forces generate moments about both  the vertical and horizontal directions. \\
			The moment about rotor shaft (vertical direction) is derived using vertical forces and moment arm $\Delta r$:
			\begin{align}
			Q   &= \frac{N}{2 \pi} \int_{0}^{2 \pi}\int_{0}^{R} \Big[\frac{\Delta D}{\Delta r}+\frac{\Delta L}{\Delta r} \frac{U_P}{U_T}\Big] r  dr \ d\Psi\\
			&= N \rho \ a \ \Bar{c}\ (\Omega R)^2 R^2 \Big[\frac{1}{8a} (1+\mu^2) \ \Bar{C_d}+ \lambda \Big(\frac{\theta_0}{6} -\frac{\theta_{tw}}{8}-\frac{\lambda}{4} \Big)\Big]    
			\end{align}
			Therefore, rotor torque coefficient $C_Q$ equals:
			\begin{eqnarray}
			\frac{C_Q}{\sigma a}   = \frac{1}{8a} (1+\mu^2) \ \Bar{C_d}+ \lambda \Big(\frac{\theta_0}{6} -\frac{\theta_{tw}}{8}-\frac{\lambda}{4}\Big)
			\end{eqnarray}
			Moments about rotor hub (horizontal direction), are rolling and pitching moments. Rolling moment is derived from vertical forces, moment arm $\Delta r$ and sine the azimuth angle $\Psi$:
			\begin{align}
			R   &= - \frac{N}{2 \pi} \int_{0}^{2 \pi}\int_{0}^{R} \frac{\Delta L}{\Delta r} r \sin{\Psi}  dr d\Psi\\
			&= -N \rho \ a \ \Bar{c} \ (\Omega R)^2 R^2 \mu \Big(\frac{\theta_0}{6} -\frac{\theta_{tw}}{8}-\frac{\lambda}{8} \Big)
			\end{align}
			The rolling coefficient is found to be:  
			\begin{eqnarray}
			\frac{C_R}{\sigma a}   = - \mu \Big(\frac{\theta_0}{6} -\frac{\theta_{tw}}{8}-\frac{\lambda}{8} \Big)
			\end{eqnarray}
			Pitching moment is derived in the same manner as rolling moment with cosine $\Psi$, instead.
			\begin{displaymath}
			P   = - \frac{N}{2 \pi} \int_{0}^{2 \pi}\int_{0}^{R} \frac{\Delta L}{\Delta r} r \cos{\Psi}  dr \ d\Psi
			\end{displaymath}
			Integration above is equal to zero, and this implies $C_P$ , the pitching coefficient, equals zero. \\
			
			Given the linear velocities and accelerations, the aerodynamic model can get the aerodynamic forces and moments on each rotor. These forces are added to the simulation in order to test the controller in a close-to-real scenario.
			\newpage
			\section{Results}
			In this section, we give results of the attitude tracking problem.  We plot the response of our control law versus the control law derived in \cite{Lee}. Here, We use the same desired signal, same moments of inertia and same initial conditions in \cite{Lee} for comparison. 
			\begin{align*}
			J &= diag(3,2,1) \\
			R(0) &= \mathbb{I} \\
			\Omega (0) &= (0,0,0)^T \\
			\phi(t) &= .999 \pi+ .5t \\
			\theta(t) &= .1 t^2 \\
			\psi(t) &= .2 t^2 -.5 t
			\end{align*}
			We use 3-2-1 Euler sequence for rotation. We choose $P$ and $F$ matrices equal to the inertia matrix.
			\begin{figure}[H]
				\centering
				\includegraphics[width = 0.75 \textwidth]{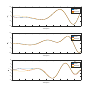}
				\label{fig:res1}
				\caption{Angular Velocity Response.}
			\end{figure}
			
			\begin{figure}[H]
				\centering
				\includegraphics[width = 0.75 \textwidth]{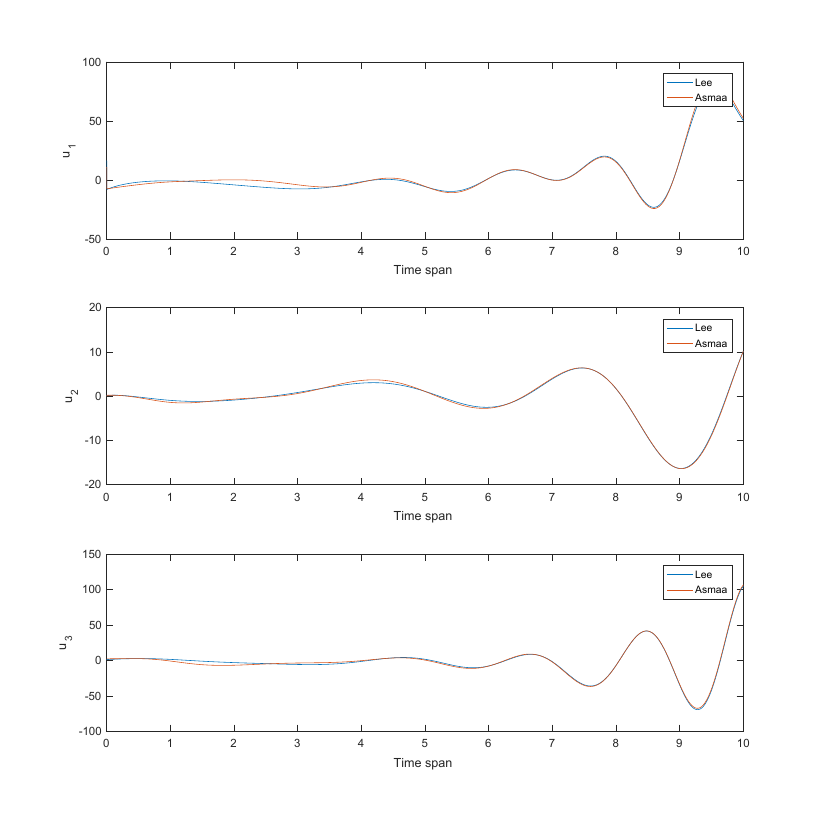}
				\label{fig:res2}
				\caption{Control Input.}
			\end{figure}
			Results show the backstepping controller gives an acceptable performance even with large initial errors in rotations close to $180^\circ$. The controller guarantees almost global asymptotic stability. Steady state errors are very small and the control effort is almost identical to the one presented in \cite{Lee}. As a consequence, we confirm that the variational integrator and the control law via backstepping are reliable. We note that a better choice for the weighing matrices  $P, F$ can even enhance the response. \\
			The simulation for the full tracking uses the following:
			\begin{align*}
			J  &= diag(0.084,0.085,0.12) \\
			m &= 4.34 \\ 
			d &= 0.315 \\
			r(0) &= (0,3,-4)^T \\
			v(0) &= (0;0,0)^T \\
			R(0) &= \mathbb{I} \\
		  \Omega (0) &= (0,0,0)^T \\
		  r_d &= 4\ (\sin(0.5t),\cos(.5t),\sin(0.5t))^T \\
		  v_d &= 2\ (\cos(.5t),-\sin(.5t),\cos(0.5t))^T \\
		  a_d &=  1\ (-\sin(.5t);-\cos(.5t);-\sin(0.5t))^T
			\end{align*}
			
			\begin{figure}[H]
					\centering
					\includegraphics[width = 0.75 \textwidth]{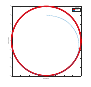}
			\end{figure}
		   \begin{figure}[H]
				\centering
				\includegraphics[width = 0.75 \textwidth]{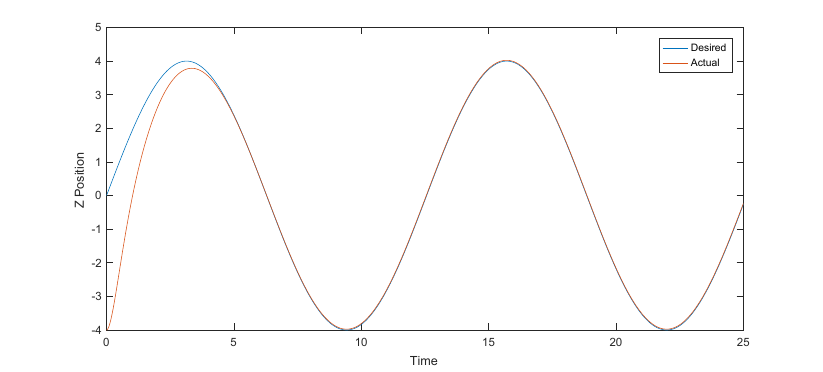}
			\end{figure}
			\begin{figure}[H]
				\centering
				\includegraphics[width = 0.75 \textwidth]{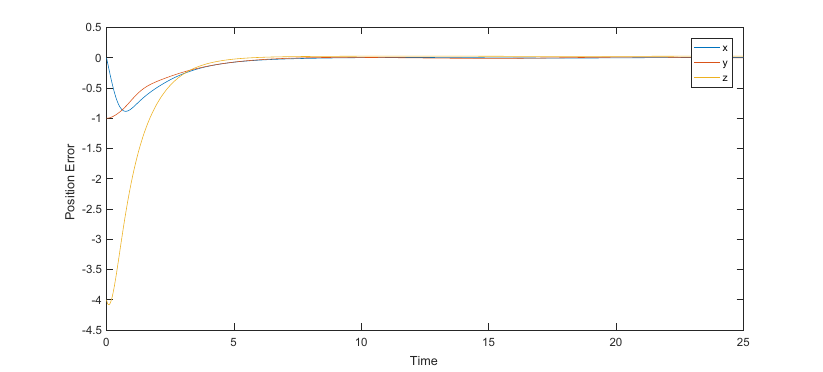}
			\end{figure}			
			
			\begin{figure}[H]
				\centering
				\includegraphics[width = 0.75 \textwidth]{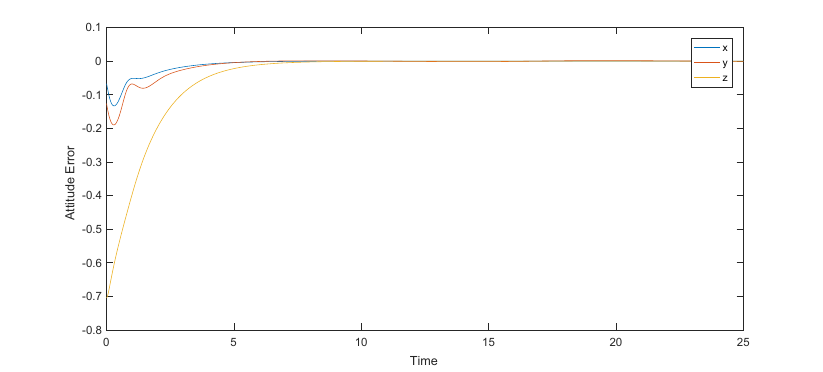}
			\end{figure}
			\begin{figure}[H]
				\centering
				\includegraphics[width = 0.75 \textwidth]{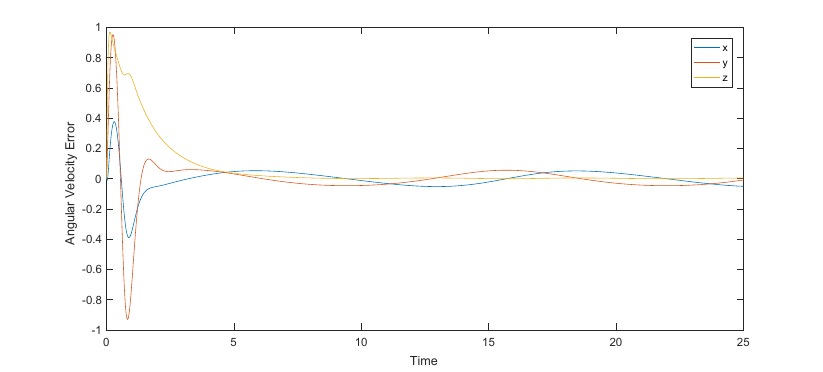}
			\end{figure}
			\begin{figure}[H]
					\centering
					\includegraphics[width = 0.75 \textwidth]{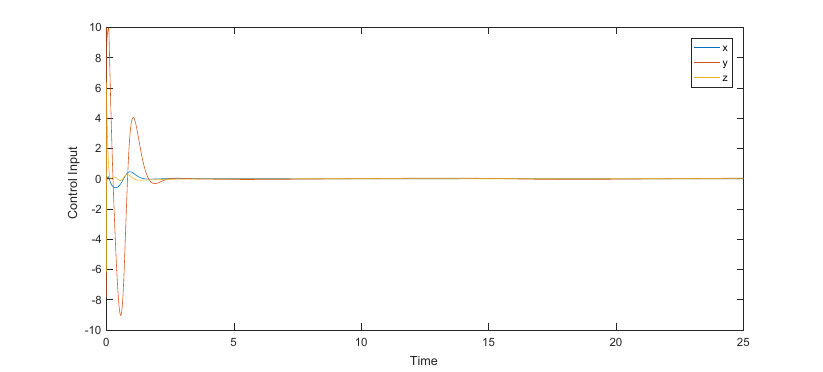}
			\end{figure}
			\begin{figure}[H]
				\centering
				\includegraphics[width = 0.75 \textwidth]{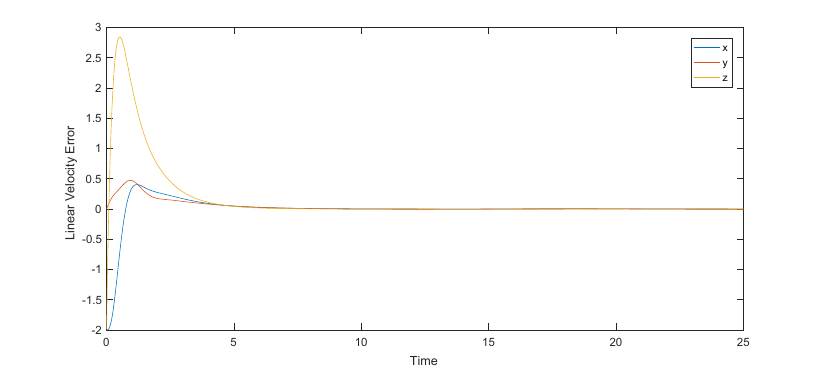}
			\end{figure}
			
			\begin{figure}[H]
				\centering
				\includegraphics[width = 0.75 \textwidth]{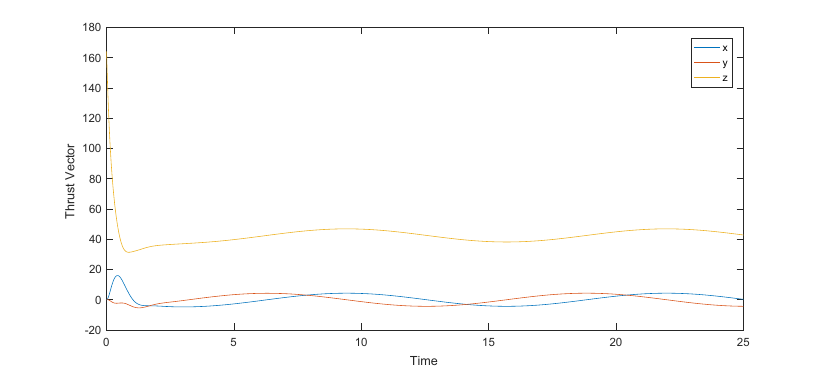}
			\end{figure}
			Results show that the quad-rotor tracks an aggressive position signal accurately with settling time around $5$ seconds. Errors in position, attitude and linear velocity converge to zero after $5$ seconds. Error in angular velocity oscillates around zero for the simulation time however it should converge as time passes. \\
			
			The effect of aerodynamics is shown in the following figures. We can see the effect obviously in tracking the z-position signal where the steady state error increases. On the other hand, response is not affected in other directions or in rotation motion.  
			\begin{figure}[H]
					\centering
					\includegraphics[width = 0.6\textwidth]{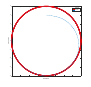}
				\end{figure}			
				
				\begin{figure}[H]
					\centering
					\includegraphics[width = 0.75 \textwidth]{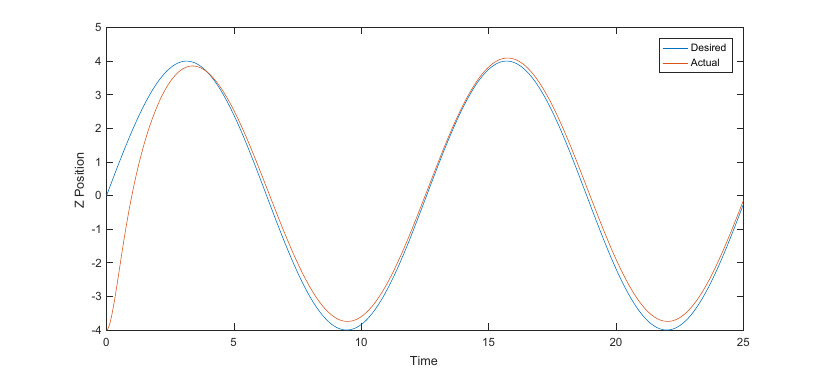}
				\end{figure}			
				\begin{figure}[H]
					\centering
					\includegraphics[width = 0.75 \textwidth]{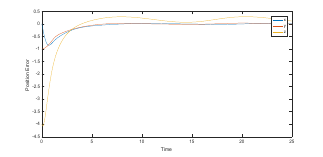}
				\end{figure}			
				
				\begin{figure}[H]
					\centering
					\includegraphics[width = 0.75 \textwidth]{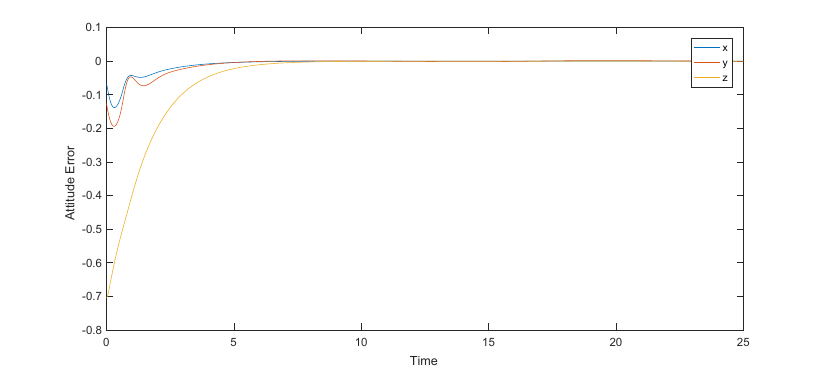}
				\end{figure}			
				
				\begin{figure}[H]
					\centering
					\includegraphics[width = 0.75 \textwidth]{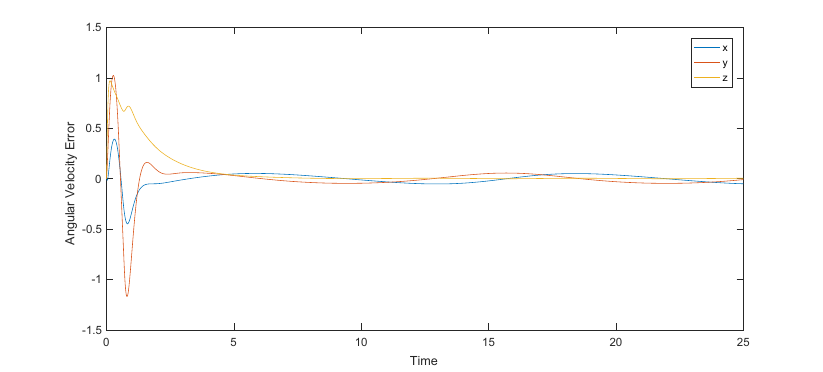}
				\end{figure}		
					
				\begin{figure}[H]
					\centering
					\includegraphics[width = 0.75 \textwidth]{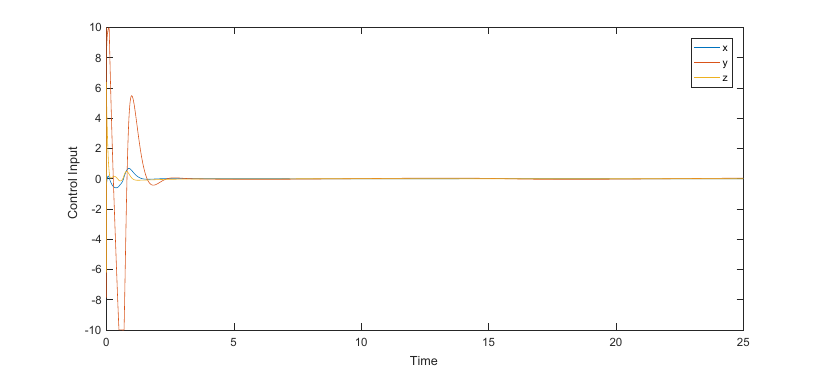}
				\end{figure}			
				
				\begin{figure}[H]
					\centering
					\includegraphics[width = 0.75 \textwidth]{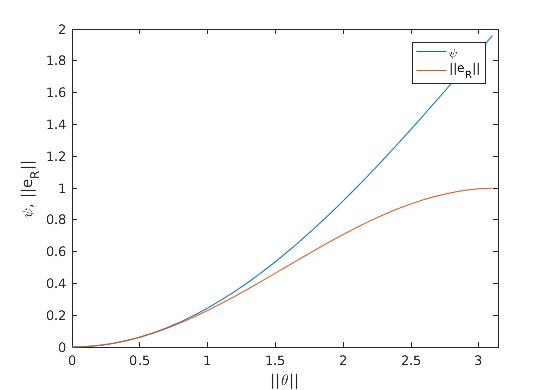}
				\end{figure}

				\begin{figure}[H]
					\centering
					\includegraphics[width = 0.75 \textwidth]{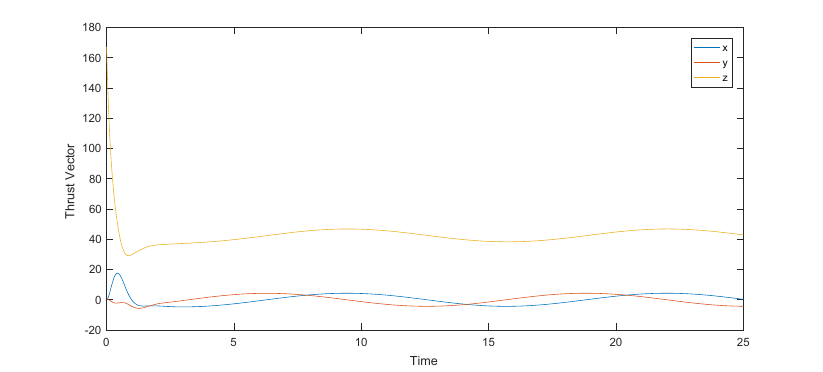}
				\end{figure}			
\section{Conclusion}
In this dissertation, we developed a 2nd-order variational integrator to rely on for rigid body simulation. With no damping in kinetic energy, we could demonstrate the ability of integration for long time. We introduced a solution for rigid body tracking with geometric backstepping technique. The proposed control laws showed very acceptable response for aggressive command maneuvers. In addition, we added the aerodynamic effects on the simulation and showed no considerable change in convergence except for the z-position where the steady state error increased.
	\bibliography {bibliography}
	\bibliographystyle{unsrt}
\end{document}